\documentclass[10pt]{amsart}

\usepackage{latexsym,amssymb}
\usepackage{amsmath,amsfonts,amsthm}
\input xypic
\xyoption{all}

\theoremstyle{plain}
\newtheorem{theorem}{Theorem}[section]
\newtheorem{proposition}[theorem]{Proposition}
\newtheorem{corollary}[theorem]{Corollary}
\newtheorem{lemma}[theorem]{Lemma}

\newtheorem{remark}[theorem]{Remark}

\begin{document}

\title[Morita equivalence of semigroups]{Morita equivalence of semigroups with local units}

\author{M.~V.~Lawson}
\address{Department of Mathematics\\
and the\\
Maxwell Institute for Mathematical Sciences\\
Heriot-Watt University\\
Riccarton\\
Edinburgh~EH14~4AS\\
Scotland}
\email{markl@ma.hw.ac.uk}

\keywords{Morita equivalence, Cauchy completions, regular semigroups.}

\subjclass[2000]{Primary: 20M10; Secondary: 20M17, 20M50.}

%\date{1.X.2009}
\begin{abstract}
We prove that two semigroups with local units are Morita equivalent if and only if they have a joint enlargement. 
This approach to Morita theory provides a natural framework for understanding McAlister's theory of the local structure of regular semigroups.
In particular, we prove that a semigroup with local units is Morita equivalent to an inverse semigroup precisely when it is a regular locally inverse semigroup.
\end{abstract}

\maketitle

%%%%%%%%%%%%%%%%%%%%%%%%%%%%%%%%%%%%%%%%%%%%%%%%%%%%%%%%%%%%%%%%%%%%%%%%%%%%%%%%%%%%%%%%%%%%%%%%%%%%%%%%%%%%%%%%%%%%%%%%%%%%%%%%%%%%%%%%%%%%%%%%%
\section{Introduction}

The Morita theory of monoids was developed independently by Banaschewski \cite{Ban} and Knauer \cite{K} and is described in \cite{KKM}.
In particular, Banaschewski showed \cite{Ban} that the generalization of this theory to semigroups cannot be accomplished by simply adjoining identities
because when this is done Morita equivalence degenerates into isomorphism.
To construct a useful Morita theory of semigroups, one has to restrict both the class of semigroups and the class of actions one considers.
The first, and decisive, step in carrying out this generalization was due to Talwar \cite{T1} who defined a Morita theory for semigroups with local units,
where a semigroup $S$ is said to have {\em local units} if for each $s \in S$ there exist idempotents $e$ and $f$ such that $es = s = sf$.
Observe that this is much weaker than the way this term is used in ring theory \cite{Abrams}.
If $S$ is a semigroup with local units then it is easy to see that $S^{2} = S$, and a semigroup with this property is said to be {\em factorizable}.
In \cite{T2,T3}, Talwar extended his theory to factorizable semigroups.
Current thinking is that factorizable semigroups form the largest class of semigroups for which a useful Morita theory can be developed.
Subsequently, only a few papers were written developing Talwar's ideas \cite{CS,N1,N2}.
Recently, however, there have been new developments.
Steinberg introduced a `strong' Morita theory for inverse semigroups \cite{S},
which turns out to be the same as the usual Morita theory of inverse semigroups, 
although in a form better adapted to inverse semigroups \cite{FL};
Laan and M\'arki \cite{LM} have been exploring Morita theory for various classes of factorizable semigroups.

In our paper, we reformulate Talwar's theory of the Morita equivalence of semigroups with local units \cite{T1} in a much more straightforward form,
and then obtain new algebraic characterizations of Morita equivalence.
As an application of our new approach, 
we show that the theory of the local structure of regular semigroups developed by McAlister \cite{M1,M1',M2,M3}
can be viewed as a contribution to the Morita theory of regular semigroups, and as a direct generalization of the pioneering paper of Rees \cite{R}.

In order to state our two main theorems, Theorems~1.1 and 1.2, we need some definitions.

We shall be dealing with actions of semigroups.\\

\noindent
{\bf Terminology }In this paper, we follow the well-established European tradition of referring to a semigroup action as an `$S$-act'
rather than as an `$S$-set'.\\

If $S$ acts on the left on the set $X$ we say that $X$ is a {\em left $S$-act}.
Left $S$-homomorphisms will be written with their arguments on the left.
Thus if $f \colon \: M \rightarrow N$ is a left $S$-homomorphism,
its value at $m$ is denoted by $(m)f$.
We denote by $S-\mbox{\bf Act}$ the category of left $S$-acts and left $S$-homomorphisms.
A left $S$-act $X$ is said to be \emph{left unitary} if and only if $SX = X$. 
If $S$ has local units and $X$ is a unitary left $S$-act, 
then it is easy to check that for each $x \in X$ there exists an idempotent $e \in S$ such that $ex = x$.
The unitary left $S$-acts with the $S$-homomorphisms between them form a full subcategory
of  $S-\mbox{\bf Act}$, which is denoted by $S-\mbox{\bf UAct}$.
If $M$ and $N$ are left $S$-acts then $\mbox{hom}_{S}(M,N)$ denotes the set
of all left $S$-homomorphisms from $M$ to $N$.
If $M$ is a right $S$-act then $\mbox{hom}_{S}(M,N)$ becomes a {\em left} $S$-act
when we define $s \cdot f$ by $(m)(s \cdot f) = (ms)f$.
In particular, $\mbox{hom}_{S}(S,M)$ is a left $S$-act.

We shall work a lot with tensor products in this paper.
Recall that two tensors $a \otimes b$ and $c \otimes d$ are equal
if there is a sequence of `moves' starting at $(a,b)$ and ending at $(c,d)$
and in each move we either move right
$(a's,b') \rightarrow (a',sb')$
or we move left
$(a',sb') \rightarrow (a's,b')$.
We can assume that left and right moves alternate by using the argument of Proposition~8.1.8 of \cite{H}
adapted to the case of semigroups with local units.

Let $X$ be a left $S$-act.
We may form the tensor product $S \otimes X$ and the action induces a map
$\mu_{X} \colon \: S \otimes X \rightarrow X$ given by $\mu_{X}(s \otimes x) = sx$.
This map is surjective if and only if $X$ is left unitary.
If it is also injective then we say that $X$ is {\em closed}.
The full subcategory of $S-\mbox{\bf Act}$ consisting of all the closed left acts is denoted by $S-\mbox{\bf FAct}$.\footnote{In Talwar's paper the `F' stands for `Fixed'
whereas for us it stands for `Ferm\'e'.} 
It is routine to check that coproducts are constructed in $S-\mbox{\bf FAct}$ in exactly the same way that they are constructed in $S-\mbox{\bf Act}$.
We define right $S$-acts dually, and we define $(S,T)$-biacts in the usual way.
A biact is {\em unitary} if it is left and right unitary.
A biact is {\em closed} if it is closed as a left and as a right act.

Let $S$ and $T$ be two semigroups  with local units. 
Then we say that $S$ and $T$ are \emph{Morita equivalent} if  the categories $S-\mbox{\bf FAct}$ and $T-\mbox{\bf FAct}$ are equivalent.
This definition is not the same as the one given by Talwar \cite{T1}, but we shall prove in Section~2 that it is equivalent to it.
This alternative definition was suggested by Neklyudova \cite{N1,N2} and is neater than the original one. 
It is easy to show that our definition coincides with the monoid one when both semigroups are monoids \cite{KKM}.

A 6-tuple $(S,T,P,Q,\langle -,- \rangle , [-,-] )$, where $S$ and $T$ are semigroups, is said to be a \emph{Morita context} if the following conditions are satisfied:
\begin{description}
\item[{\rm (M1)}] $P$ is an $(S,T)$-biact, and $Q$ is an $(T,S)$-biact.

\item[{\rm (M2)}] $\langle -,-\rangle \colon \: P \otimes Q \rightarrow S$ is an  $(S,S)$-homomorphism
and  $[-,-] \colon \: Q \otimes P \rightarrow T$ is a  $(T,T)$-homomorphism.

\item[{\rm (M3)}] For all $p,p' \in P, q,q' \in Q$ the following two conditions are satisfied:
\begin{description}
\item[{\rm (i)}] $\langle p,q \rangle p' = p[q,p']$.
\item[{\rm (ii)}] $q\langle p,q'\rangle = [q,p]q'$.
\end{description}
 \end{description}
We say that a Morita context $(S,T,P,Q,\langle -,-\rangle,[-,-])$ is \emph{unitary} if and only if 
$S$ and $T$ are semigroups with local units,
$P$ and $Q$ are closed as left acts, 
and the biacts $P$ and $Q$ are unitary.

A semigroup $S$ is {\em regular} if for each $s \in S$ there exists $t \in S$ such that $s = sts$ and $t = tst$.
The element $t$ is called an {\em inverse} of $S$.
The set of inverses of $s$ is denoted by $V(s)$.
If each element has a unique inverse then the semigroup is said to be {\em inverse}.
For undefined terms from regular semigroup theory see \cite{H}.

Let $S$ be a subsemigroup of the semigroup $T$.
Then $T$ is said to be an {\em enlargement} of $S$ if $S = STS$ and $T = TST$.
Let $S$, $T$ and $R$ be semigroups with local units.
We shall say that $R$ is a {\em joint enlargement} of $S$ and $T$
if it is an enlargement of subsemigroups $S'$ and $T'$ which are isomorphic to $S$ and $T$ respectively.
The theory of enlargements was introduced in \cite{L1} and developed in \cite{L3}.
Steinberg's paper \cite{S} was explicitly motivated by enlargements.

Categories will be used both as structures on a par with monoids as well as the more usual categories of structures;
it will be clear from the context which perspective is intended.
Furthermore, definitions from semigroup theory can be extended in the obvious way to categories.
If $S$ is a semigroup then
$$C(S) = \{(e,s,f) \in E(S) \times S \times E(S) \colon \: esf = s \}$$
is a category called the {\em Cauchy completion} of $S$.
We shall build semigroups from (small) categories using the following technique.
A category $C$ is said to be {\em strongly connected} if for each pair of identities $e$ and $f$ there is an arrow from $e$ to $f$.
Let $C$ be a strongly connected category.
A {\em consolidation} for $C$ is a function $p \colon \: C_{o} \times C_{o} \rightarrow C$,
$p(e,f) = p_{e,f}$, where $p_{e,f}$ is an arrow from $f$ to $e$ and $p_{e,e} = e$.
Given a category $C$ equipped with a consolidation $p$ we can define a binary operation $\circ$ on $C$ by
$x \circ y = xp_{e,f}y$ where $x$ has domain $e$ and $y$ has codomain $f$.
It is easily checked that this converts $C$ into a semigroup.
We denote this semigroup by $C^{p}$.
If we omit $\circ$ then the product is in the category.

Let $S$ and $T$ be semigroups with local units.
A homomorphism $\theta \colon \: S \rightarrow T$ is said to be a {\em local isomorphism} if the following conditions are satisfied:
\begin{description}

\item[{\rm (LI1)}] The function $\theta$ restricted to $eSf$
induces an isomorphism with $\theta (e)T \theta (f)$ for all idempotents $e$ and $f$ in $S$.

\item[{\rm (LI2)}] {\em Idempotents lift along} $\theta$ meaning that if $e'$ is an idempotent in 
the image of $\theta$ then there is an idempotent $e$ in $S$ such that $\theta (e) = e'$.

\item[{\rm (LI3)}] For each idempotent $e \in T$ there exists an idempotent $f \in T$ in the image of $\theta$ such that $e \, \mathcal{D} \, f$.

\end{description}
This is a generalization of the classical definition of a local isomorphism between regular semigroups \cite{M1,M1'}
and has its origins in \cite{MS} and \cite{L3} as well as topos theory.
When $S$ is regular, 
surjective local isomorphisms are precisely the surjective homomorphisms that are injective when restricted to each local submonoid \cite{L3}.

We shall prove two main theorems in this paper.
The first describes different characterizations of Morita equivalence.

\begin{theorem} Let $S$ and $T$ be semigroups with local units.
Then the following are equivalent.
\begin{enumerate}

\item $S$ and $T$ are Morita equivalent.

\item The categories $C(S)$ and $C(T)$ are equivalent

\item $S$ and $T$ have a joint enlargement which can be chosen to be regular if $S$ and $T$ are both regular.

\item There is a unitary Morita context $(S,T,P,Q,\langle -,- \rangle,[-,-])$ with surjective mappings.

\end{enumerate}
\end{theorem}

The second describes a practical starting point for trying to show that two semigroups are Morita equivalent.
It adapts to our setting the heuristic described by McAlister in \cite{Mc}.

\begin{theorem} Let $S$ and $T$ be semigroups with local units.
Then $S$ and $T$ are Morita equivalent if and only if 
there is a consolidation $q$ on $C(S)$ and a local isomorphism $\psi \colon \: C(S)^{q} \rightarrow T$.
\end{theorem} 

In Section~2, we shall reconcile our approach with Talwar's and, apart from Proposition~2.4, 
we do not essentially use the results of this section later.
In Sections~3 and 4, we prove Theorems~1.1 and 1.2, respectively.
In Section~5, we apply our results to McAlister's theory of the local structure of regular semigroups.

One question, raised by the referee, which we do not solve here, is the following: if $S$ and $T$ are Morita equivalent semigroups
what can we say about their semigroup rings $RS$ and $RT$ where $R$ is a commutative ring with identity?
The case where $S$ and $T$ have commuting idempotents is easy to handle:
the semigroup rings are Morita equivalent.
The proof is the same as that of Theorem~4.13 of \cite{S};
one observes that $RS$ and $RT$ are rings with local units in the sense of Abrams \cite{Abrams} and then 
constructs from the semigroup Morita context guaranteed by Theorem~1.1(4) a ring Morita context.\\

\noindent
{\bf Acknowledgements }I am grateful to Jonathon Funk and Benjamin Steinberg for numerous email exchanges on the subject of this paper.
In particular, the idea that the main implication (2)$\Rightarrow$(3) of Theorem~1.1 should be true arose out of discussions with them.
The definition of `local isomorphism' generalizes a definition made by Jonathon Funk in the inverse case and originates in topos theory. 
I am also grateful to Valdis Laan and Laszlo M\'arki for helpful discussions, for generously sharing their work-in-progress \cite{LM},
and for spotting some errors and typos. 
The referee helpfully suggested a number of improvements including aspects of the proofs of Theorem~2.1 and Proposition~5.1.
Finally, my overriding thanks must go to Don McAlister whose research has been an inspiration for my own.

%%%%%%%%%%%%%%%%%%%%%%%%%%%%%%%%%%%%%%%%%%%%%%%%%%%%%%%%%%%%%%%%%%%%%%%%%%%%%%%%%%%%%%%%%%%%%%%%%%%%%%%%%%%%%%%%%%%%%%%%%%%%%%%%%%%%%%%%%%%%%%%%%%%%%
\section{The category of closed left acts} 

In this section, we explain the connection between our approach to Morita theory and the one pioneered by Talwar.
We begin by discussing a couple of minor problems in Talwar's account.
In his paper, Talwar \cite{T1} defines the class of unitary left $S$-acts $X$ to be considered by requiring that
the evaluation map $S \otimes \mbox{hom}(S,X) \rightarrow X$ be an isomorphism.
Towards the end of the paper he proves that such acts are precisely those for which there is {\em an} isomorphism $S \otimes X \cong X$.
However, when showing that two categories are equivalent one needs {\em natural} isomorphisms.
For this reason, one should work with the natural isomorphisms $\mu_{X} \colon \: S \otimes X \rightarrow X$ defined in Section~1.
This agrees with what is done in ring theory \cite{M}.
In Proposition~2.2, we prove that nevertheless our definition of Morita equivalence coincides with Talwar's.
Another minor problem with Talwar's paper is that he assumes epimorphisms are surjective.
This is easily rectified in our Proposition~2.4.

There are two categories of interest to us:
$$S-\mbox{\bf FAct} \subseteq S-\mbox{\bf UAct}$$
where the inclusion is as full subcategories.
We shall define two endofunctors of $S-\mbox{\bf UAct}$ which we will use to better understand $S-\mbox{\bf FAct}$.

\begin{itemize}

\item The functor $S \otimes \colon \:  S-\mbox{\bf Act} \rightarrow  S-\mbox{\bf UAct}$ is defined in the usual way.
Let $s \otimes m \in S \otimes M$, and let $e$ be an idempotent such that $es = s$.
Then $e(s \otimes m) = es \otimes m = s \otimes m$.
Thus  $S \otimes M$ is always unitary.
We have already defined $\mu_{M}$ in Section~1.
These form the components of a natural transformation $\mu$ from the functor $S \otimes -$ to the identity functor on the category $S-\mbox{\bf UAct}$.

\item The functor $S \colon \:  S-\mbox{\bf Act} \rightarrow  S-\mbox{\bf UAct}$
is defined as follows.
Let $M$ be any $S$-act. Then the set  $$SM = \{ sm \colon s \in S, m \in M \}$$ is a
unitary subact of $M$. In addition, if $N$ is a unitary subact of $M$ then $N \subseteq SM$.
Thus $SM$ is the largest unitary subact of $M$.
If $f \colon \: M \rightarrow N$ is any left $S$-homomorphism
then it restricts to a left $S$-homomorphism $Sf \colon \: SM \rightarrow SN$ given by $(m)(Sf) = (m)f$.
It is clearly right adjoint to the forgetful functor from  $S-\mbox{\bf UAct}$ to $S-\mbox{\bf Act}$.

\item The functor $\mbox{\rm hom}_{S}(S,-) \colon \:  S-\mbox{\bf Act} \rightarrow  S-\mbox{\bf Act}$ is defined in the usual way.
If $f \colon \: M \rightarrow N$ is a left $S$-homomorphism then we denote
$$\mbox{\rm hom}_{S}(f) \colon \: \mbox{\rm hom}_{S}(S,M) \rightarrow  \mbox{\rm hom}_{S}(S,N),$$
which maps $\alpha$ to $\alpha f$, by $f^{\ast}$.
For each $m \in M$, define 
$$\rho_{m} \colon \: S \rightarrow M$$ 
by $(s)\rho_{m} = sm$. 
This is a left $S$-homomorphism.
The function 
$$\rho_{M} \colon \: M \rightarrow \mbox{\rm hom}_{S}(S,M)$$
defined by $m \mapsto \rho_{m}$ is a left $S$-homomorphism
and forms the components of a natural transformation $\rho$ 
from the identity functor on $S-\mbox{\bf Act}$ to the functor $\mbox{\rm hom}_{S}(S,-)$.

\item The functor $S\mbox{\rm hom}_{S}(S,-) \colon \:  S-\mbox{\bf Act} \rightarrow  S-\mbox{\bf UAct}$ 
combines the above two functors and is an endofunctor of  $S-\mbox{\bf UAct}$.
If $M$ is unitary then the image of $\rho_{M} \colon \: M \rightarrow \mbox{\rm hom}_{S}(S,M)$
is a unitary subact and so we can regard it as a left $S$-homomorphism
$\rho_{M} \colon \: M \rightarrow S\mbox{\rm hom}_{S}(S,M)$ and this forms the components of
a natural transformation from the identity functor on $S-\mbox{\bf UAct}$ to the functor $S\mbox{\rm hom}_{S}(S,-)$.
\end{itemize}

The two endofunctors of $S-\mbox{\bf UAct}$ defined above are related by the following theorem
which implies that $S \otimes -$ and $S\mbox{hom}_{S}(S,-)$ form a {\em Galois adjunction}
on the category $S-\mbox{\bf UAct}$; see Chapter~19, Exercise~19D of \cite{AHS}.
The isomorphism of (5) below is proved as Lemma~4.8 of \cite{T1}.

\begin{theorem} Let $S$ be a semigroup with local units.
Then on the category $S-\mbox{\bf UAct}$, we have the following.
\begin{enumerate}

\item The functor $S \otimes -$ is left adjoint to the functor $S\mbox{\rm hom}_{S}(S,-)$.

\item The unit of the adjunction is the function
$$\eta_{M} \colon \: M \rightarrow S\mbox{\rm hom}_{S}(S,S \otimes M)$$
given by $m \mapsto - \otimes m$.

\item The counit of the adjunction is the function 
$$\varepsilon_{M} \colon \: S \otimes S\mbox{\rm hom}_{S}(S,M) \rightarrow M$$
given by $s \otimes f \mapsto (s)f$.

\item $S \otimes S\mbox{\rm hom}_{S}(S,M) =  S \otimes \mbox{\rm hom}_{S}(S,M)$. 

\item The function $S \otimes S\mbox{\rm hom}_{S}(S,S \otimes M)   \stackrel{\varepsilon_{S \otimes M}}{\rightarrow} S \otimes M$ 
is an isomorphism with inverse $1 \otimes \eta_{M}$.

\end{enumerate}
\end{theorem}
\begin{proof} The forgetful functor $U \colon  \:  S-\mbox{\bf UAct} \rightarrow  S-\mbox{\bf Act}$ 
has right adjoint the functor $S \colon \:  S-\mbox{\bf Act} \rightarrow  S-\mbox{\bf UAct}$,
and the functor $S \otimes \colon \:  S-\mbox{\bf Act} \rightarrow  S-\mbox{\bf UAct}$ 
has, as usual, the right adjoint $\mbox{\rm hom}_{S}(S,-) \colon \:  S-\mbox{\bf UAct} \rightarrow  S-\mbox{\bf Act}$.
However adjunctions compose \cite{Mac}.
This proves (1).
The proofs of (2) and (3) are now routine.

(4) To show that $S \otimes S\mbox{\rm hom}_{S}(S,M) =  S \otimes \mbox{\rm hom}_{S}(S,M)$. 
Observe that $S \otimes S\mbox{\rm hom}_{S}(S,M) \subseteq  S \otimes \mbox{\rm hom}_{S}(S,M)$. 
To prove the reverse inclusion, let $s \otimes f \in S \otimes \mbox{\rm hom}_{S}(S,M)$. 
Let $e^{2} = e$ be such that $se = s$.
Then $s \otimes f = s \otimes e \cdot f$.
But $e \cdot f \in  S\mbox{\rm hom}_{S}(S,M)$.

(5) From Theorem~IV.1 of \cite{Mac}, there is a left $S$-homomorphism given by
$$1 \otimes \eta_{M} \colon \: S \otimes M \rightarrow S \otimes S \mbox{hom}_{S}(S,S \otimes M).$$
The effect of this function is
$s \otimes m \mapsto s \otimes \rho_{e \otimes m}$ 
where $e$ is any idempotent such that $em = m$.
We also have a left $S$-homomorphism going the other way
$$\varepsilon_{S \otimes M} \colon \: S \otimes S \mbox{hom}_{S}(S,S\otimes M) \rightarrow S \otimes M,$$
given by $s \otimes f \colon \mapsto (s)f$.
It follows from the general theory of adjunctions, 
and can easily be directly verified,
that $(1 \otimes \eta_{M})\varepsilon_{S \otimes M}$ is the identity.
Let 
$s \otimes f \in S \otimes \mbox{\rm hom}_{S}(S,S \otimes M)$.
We calculate
$$(s \otimes f)\varepsilon_{S \otimes M}(1 \otimes \eta_{M}).$$
Now $(s)f \in S \otimes M$ which means that $(s)f = t \otimes m$ for some $t \in S$ and $m \in M$. 
Let $e^2 = e \in S$ such that $es = s$, and let $f^{2} = f$ be such that $fm = m$. 
Then 
$$t \otimes m = (s)f = (es)f = e(s)f = e(t \otimes m).$$
Now
\begin{eqnarray*}
(t \otimes m)(1 \otimes \eta_{M})  &=& (et \otimes m)(1 \otimes \eta_{M})\\
                                   &=& et \otimes \rho_{f \otimes m}\\
                                  &=& e \otimes t \cdot \rho_{f \otimes m}\\ 
                                   &=& e \otimes \rho_{tf \otimes m}\\ 
                                  &=& e \otimes \rho_{t \otimes m}\\
                                  &=& e \otimes \rho_{(s)f}\\
                                  &=& e \otimes s \cdot f\\
                                  &=& s \otimes f
\end{eqnarray*}
as required.
\end{proof}

%%%%%%%%%%%%%%%%%%%%%%%%%%%%%%%%%%%%%%%%%%%%%%%%%%%%%%%%%%%%%%%%%%%%%%%%%%%%%%%%%%%%%%%%%%%%%%%%%%%%%%%%%%%%%%%%%%%%%%%%%
\begin{lemma}
There is a  left $S$-isomorphism
$$1 \otimes \rho_{M} \colon \: S \otimes M \rightarrow S \otimes S \mbox{\rm hom}_{S}(S,M)$$
defined by $s \otimes m \mapsto s \otimes \rho_{m}$ which is natural in $M$.
\end{lemma}
\begin{proof}
We show first that this is a well-defined function.
Map the ordered pair $(s,m)$ to $s \otimes \rho_{m}$.
Thus $(st,m)$ maps to $st \otimes \rho_{m}$ whereas $(s,tm)$ maps to $s \otimes \rho_{tm}$.
However
$$(s')\rho_{tm} = (s')(tm) = (s't)m = (s')(t \cdot \rho_{m}).$$
Thus $\rho_{tm} = t \cdot \rho_{m}$.
But then we have that $s \otimes \rho_{tm} = s \otimes t \cdot \rho_{m} = st \otimes \rho_{m}$.
It follows that the function $s \otimes m \mapsto s \otimes \rho_{m}$ is well-defined.
It is therefore clear that we have defined a left $S$-homomorphism.

We now define a function going in the other direction.
Map the ordered pair $(s,\alpha)$ to $e \otimes (s)\alpha$ where $e$ is any idempotent such that $es = s$.
We prove first that this is a well-defined function; that is, independent of the choice of idempotent $e$.
Let $si = s$ where $i$ is an idempotent.
Then 
$$e \otimes (s)\alpha = e \otimes (si)\alpha = e \otimes s(i)\alpha = es \otimes (i)\alpha = s \otimes (i)\alpha.$$
Now let $fs = s$ where $f$ is any idempotent.
Then
$$s \otimes (i)\alpha = fs \otimes (i)\alpha = f \otimes s(i)\alpha = f \otimes (si)\alpha = f \otimes (s)\alpha.$$
In this case, it is easy to check that $(st,\alpha)$ and $(s,t \cdot \alpha)$ map to the same element.
It follows that we have a well-defined funcion
$$S \otimes S \mbox{\rm hom}_{S}(S,M) \rightarrow S \otimes M$$ 
given by $s \otimes \alpha \mapsto e \otimes (s)\alpha$ where $es = s$ is any idempotent.
It can now be checked that this map is a left $S$-homomorphism.

We now show that these two left $S$-homomorphisms are mutually inverse.
Let $s \otimes m \in S \otimes M$.
Then this maps to $s \otimes \rho_{m}$ which in turn maps to $e \otimes (s)\rho_{m}$ where $es = s$.
But $e \otimes (s)\rho_{m} = e \otimes sm = es \otimes m  = s \otimes m$.

Let $s \otimes \alpha \in S \otimes \mbox{\rm hom}_{S}(S,M)$.
Then this maps to $e \otimes (s)\alpha$ where $es = s$ which in turn maps to $e \otimes \rho_{(s)\alpha}$.
It is easy to check that $\rho_{(s)\alpha} = s \cdot \alpha$.
Thus 
$$e \otimes \rho_{(s)\alpha} = e \otimes s \cdot \alpha = es \otimes \alpha = s \otimes \alpha.$$
\end{proof}

%%%%%%%%%%%%%%%%%%%%%%%%%%%%%%%%%%%%%%%%%%%%%%%%%%%%%%%%%%%%%%%%%%%%%%%%%%%%%%%%%%%%%%%%%%%%%%%%%%%%%%%%%%%%%
It follows from Theorem~2.1 and Lemma~2.2 that 
$S-\mbox{\bf FAct}$ is a full coreflective subcategory of $S-\mbox{\bf UAct}$:
the coreflection of the unitary left $S$-act $X$ is $S \otimes X$ and there is an epimorphism
$\mu_{X} \colon \: S \otimes X \rightarrow X$.
We are now able to prove that our definition of Morita equivalence is the same as Talwar's.

\begin{proposition} Let $S$ be a semigroup with local units, and let $M$ be a unitary left $S$-act.
Then the following are equivalent.
\begin{enumerate}

\item $M$ is closed.

\item $S \otimes M \cong M$; that is, $S \otimes M$ is isomorphic to $M$ for some isomorphism.

\item The map 
$$\varepsilon_{M} \colon \: S \otimes S \mbox{\rm hom}(S,M) \rightarrow M$$
is an isomorphism.

\end{enumerate}
\end{proposition}
\begin{proof}

(1)$\Rightarrow$(2). This is immediate.

(2)$\Rightarrow$(3). Let $f \colon \: S \otimes M \rightarrow M$ be a left $S$-isomorphism.
The following diagram commutes
$$\xymatrix@!{
&S \otimes S \mbox{hom}_{S}(S,S\otimes M) \ar[r]^{\varepsilon_{S \otimes M}} \ar[d]_{1 \otimes f^{\ast}}
&S \otimes M \ar[d]^{f}
\\
&S \otimes S \mbox{hom}_{S}(S,M) \ar[r]^{\varepsilon_{M}}
&M
}$$
By assumption, $f$ is an isomorphism and so $1 \otimes f^{\ast}$ is an isomorphism.
By Theorem~2.1(5),
$\varepsilon_{S \otimes M} \colon \: S \otimes \mbox{hom}_{S}(S,S \otimes M) \rightarrow S \otimes M$
is an isomorphism.
Thus $\varepsilon_{M}$ is an isomorphism, as required.

(3)$\Rightarrow$(1). Suppose that $\varepsilon_{M}$ is an isomorphism.
The following diagram commutes
$$\xymatrix@!{
&S \otimes M \ar[r]^{1 \otimes \rho_{M}} \ar[d]_{\mu_{M}} 
&S \otimes S \mbox{hom}_{S}(S,M) \ar[dl]^{\varepsilon_{M}}
\\
&M 
&
}$$
By Lemma~2.2, $1 \otimes \rho_{M}$ is an isomorphism.
It follows that $\mu_{M}$ is an isomorphism, and so $M$ is closed.
\end{proof}

%%%%%%%%%%%%%%%%%%%%%%%%%%%%%%%%%%%%%%%%%%%%%%%%%%%%%%%%%%%%%%%%%%%%%%%%%%%%%%%%%%%%%%%%%%%%%%%%%%%%%%%%%%%%%%%%%%%%%%%%%%%%%%
We conclude this section by proving that epimorphisms are always surjective in the category of closed left acts.
Let $M$ be a unitary left $S$-act.
An equivalence relation $\sim$ on $M$ is said to be a {\em left congruence} if
$m \sim n$ implies that $sm \sim sn$ for all $s \in S$.
Denote the $\sim$-equivalence class containing $m$ by $[m]$.
Then $M/\sim$ is also unitary left $S$-act.
The intersection of left congruences on $M$ is again a left congruence,
so we can talk about the left congruence generated by a relation.  
The proof of the following is adapted from \cite{Ban}.

\begin{proposition} In the category $S-\mbox{\bf FAct}$, all epimorphisms are surjections.
\end{proposition}
\begin{proof}
Let $f \colon \: M \rightarrow N$ be an epimorphism in the category $S-\mbox{\bf FAct}$.
Then we have the following diagram
$$\diagram
&M \ar[r]^{f} 
&N
\\
&S \otimes M \ar[r]^{1 \otimes f} \ar[u]^{\mu_{M}} 
&S \otimes N \ar[u]_{\mu_{N}}
\enddiagram$$
which commutes.
Since $\mu_{M}$ and $\mu_{N}$ are isomorphisms it follows that $1 \otimes f$ is an epimorphism.
Let the image of $f$ be the left $S$-subact $N'$ of $N$.
We shall suppose that $N' \neq N$ and derive a contradiction
from which it will follow that $f$ is surjective.
Form the coproduct $N \sqcup N = N \times \{1\} \cup N \times \{ 2\}$.
The elements of this coproduct are of the form $(n,i)$ where $i = 1,2$ and
the left $S$-action is given by $s(n,i) = (sn,i)$.
Define the relation $\sim$ on $N \sqcup N$ by
$(x,i) \sim (x',j)$ iff either $(x,i) = (x',j)$ or $i \neq j$ and $x = x' \in N'$.
This is an equivalence relation on $N \sqcup N$ and a left congruence.
The $\sim$-equivalence class containing $(x,i)$ is denoted by $[(x,i)]$.
We denote the set of $\sim$-equivalence classes by $N_{f}$.
There are two left $S$-homomorphisms
$j_{1},j_{2} \colon \:N \rightarrow N_{f}$ given by
$(n)j_{1} = [(n,1)]$ and $(n)j_{2} = [(n,2)]$.
Observe that
$$fj_{1} = fj_{2}$$
but 
$$(n)j_{1} \neq (n)j_{2}$$
for any $n \in N \setminus N'$.
Observe that
$$(1 \otimes f)(1 \otimes j_{1}) = (1 \otimes f)(1 \otimes j_{2})$$
and that $1 \otimes f$ is an epimorphism;
all these maps are in the category $S-\mbox{\bf FAct}$.
It therefore only remains to prove that
$1 \otimes j_{1} \neq 1 \otimes j_{2}$
to derive our contradiction.
Let $n \in N \setminus N'$,
and let $e$ be an idempotent in $S$ such that $en = n$.
Then $e \otimes n \in S \otimes N$.
Suppose that
$(e \otimes n)(1 \otimes j_{1}) = (e \otimes n)(1 \otimes j_{2})$.
Then $e \otimes (n)j_{1} = e \otimes (n)j_{2}$.
Applying the map $\mu_{N_{f}}$ we get $e(n)j_{1} = e(n)j_{2}$
and so $(n)j_{1} = (n)j_{2}$, which is a contradiction.
It follows that $N \setminus N'$ is empty and so $f$ is a surjection, as required.
\end{proof}

%%%%%%%%%%%%%%%%%%%%%%%%%%%%%%%%%%%%%%%%%%%%%%%%%%%%%%%%%%%%%%%%%%%%%%%%%%%%%%%%%%%%%%%%%%%%%%%%%%%%%%%%%%%%%%%%%%%%%%%%%%%%%%%%%%%%%%%%%%%%%%%%%%%%%%%%5
\section{Proof of Theorem~1.1}

In this section, we shall prove each of the implications in the statement of Theorem~1.1.

%%%%%%%%%%%%%%%%%%%%%%%%%%%%%%%%%%%%%%%%%%%%%%%%%%%%%%%%%%%%%%%%%%%%%%%%%%%%%%%%%%%%%%%%%%%%%%%%%%%%%%%%%%%%%%%%%%%%%%%%%%%%%%%%%%%%%%%%%%%%%%%%%%%%%
\subsection{From Morita equivalence to Cauchy completions}

The main result of this subsection, Theorem~3.4, was known to Talwar \cite{T2}, and is included here for the sake of completeness.
We say that a closed left $S$-act $M$ is {\em indecomposable} if
$M$ is not isomorphic to any coproduct $N \sqcup N'$ where $N$ and $N'$ are non-empty
closed left $S$-acts.

\begin{lemma} Let $S$ be a semigroup with local units.
\begin{enumerate}

\item For each idempotent $e \in S$ the left $S$-act $Se$ is closed.

\item In the category $S-\mbox{\bf FAct}$ the closed left $S$-acts of the form $Se$, where $e$ is an idempotent, are indecomposable and projective.

\end{enumerate}
\end{lemma}
\begin{proof}
(1) It is clear that $Se$ is a unitary left $S$-act.
We shall show that the function $\mu \colon \: S \otimes Se \rightarrow Se$ is injective.
Suppose that $(s \otimes a)\mu = (t \otimes b)\mu$ where $s,t \in S$ and $a,b \in Se$.
Then $sa = tb$.
Let $fs = s$. 
Then $ftb = tb$.
Thus
$$s \otimes a = fs \otimes a = f \otimes sa = f \otimes tb$$ 
but $tb = (tb)e$ and so
$$f \otimes tb = f \otimes (tb)e = f(tb) \otimes e = tb \otimes e = t \otimes be = t \otimes b,$$
as required.

(2) Suppose that $M$ and $N$ are two, non-empty, closed left $S$-subacts of $Se$ such that
$Se = M \cup N$ and $M \cap N = \emptyset$.
Now $e \in Se$ and so $e \in M$ or $e \in N$.
Without loss of generality, assume the former.
Then $e \in M$ implies that $Se \subseteq M$.
Thus $N \subseteq M \cup N = Se \subseteq M$ and so $N \subseteq M$, which is a contradiction.
Thus $Se$ is indecomposable.

To prove that $Se$ is projective.
Let $f \colon \: M \rightarrow N$ be an epimorphism and let
$g \colon \: Se \rightarrow N$ be arbitrary.
By Proposition~2.4, $f$ is surjective and so there is $a \in M$ such that $(a)f = (e)g$.
Define $h \colon \: Se \rightarrow M$ by $(se)h = sea$.
Then $h$ is a left $S$-homomorphism.
Now $(se)(hf) = (sea)f = se(a)f = se(e)g = (se)g$.
Thus $hf = g$, as required.
\end{proof}

In the next lemma, we assemble some results on projectives in the category $S-\mbox{\bf FAct}$.

\begin{lemma} In the category $S-\mbox{\bf FAct}$ the following hold.
\begin{enumerate}

\item Every coproduct of projectives is projective.

\item The category has enough projectives.

\item Let $\theta \colon \: P \rightarrow P'$ and $\theta' \colon \:P' \rightarrow P$ be such that
$\theta \theta' = 1_{P}$. Then if $P'$ is projective so is $P$.

\item $P$ is projective if and only if every epimorphism $M \rightarrow P$ has a left inverse.

\end{enumerate}
\end{lemma}
\begin{proof}
(1) An easy deduction from Proposition~14.3 of \cite{M}.

(2) Let $M$ be an arbitrary closed left $S$-act.
For each $m \in M$ choose an idempotent $e_{m}$ such that $e_{m}m = m$.
Form the coproduct $\coprod_{m \in M}Se_{m}$.
This is projective and unitary by Lemma~3.1 and (1) above, and closed because the coproduct of closed acts is closed.
Define $\pi \colon \: \coprod_{m \in M}Se_{m} \rightarrow M$ by
$(se_{m})\pi = se_{m}m = sm$.
This is a surjective left $S$-homomorphism.
For convenience, we shall call the map $\pi$ defined above the {\em canonical covering} of $M$.

(3) Proposition~14.1 of \cite{M}.

(4 ) An easy deduction from Proposition~14.2 of \cite{M} and (2) above.
\end{proof}

\begin{proposition} In the category $S-\mbox{\bf FAct}$,
$P$ is indecomposable and projective if and only if $P$ is isomorphic to $Se$ for some idempotent $e$.
\end{proposition}
\begin{proof}
Lemma~3.1 proves one direction so we need only prove the converse.
Let $P$ be indecomposable and projective.
By Lemma~3.2(2), there is the canonical covering
$\pi \colon \: \coprod_{p \in P}Se_{p} \rightarrow P$.
By Lemma~3.2(4), there is an injective left $S$-homomorphism
$\sigma \colon \: P \rightarrow \coprod_{p \in P}Se_{p}$ such that $\sigma \pi = 1_{P}$.
Now $P$ is indecomposable by Lemma~3.1, 
and $(P)\sigma$ is isomorphic to $P$ and so also indecomposable.
But $(P)\sigma $ is a subact of $\coprod_{p \in P}Se_{p}$.
It follows that $(P)\sigma \subseteq Se_{p}$ for some $p \in P$.
Thus $\sigma \colon \: P \rightarrow Se_{p}$ defines an injective left $S$-homomorphism.
But using the fact that $\sigma \pi = 1_{P}$ we find that $P = (Se_{p})\pi$.
Now $Se_{p}$ is a cyclic left $S$-act and so $P$ is a cyclic left $S$-act.
We may therefore assume that $P$ is a projective cyclic left $S$-act where $P = Sx$ for some $x \in P$.
Since $P$ is closed it is, in particular, unitary and so there is an idempotent $e \in S$ such that $ex = x$.
Define $\phi \colon \: Se \rightarrow P$ by $(s)\phi = sx$.
Then $\phi$ is a surjection.
But $Sx$ is projective and so there exists a map $\psi \colon \: P \rightarrow Se$ such that $\psi \phi = 1_{P}$.
We therefore have an injective map $\psi \colon \: P = Sx \rightarrow Se$.
Put $f = (x)\psi$.
Then $f = (x)\psi = (ex)\psi = e(x)\psi = ef$, and since $x \in Sf$ we have that $fe = f$.
Observe that $f^{2} = fefe = fe = f$ and so $f$ is an idempotent and $f \leq e$.
It follows that $\psi$ induces an isomorphism between $P = Sx$ and $Se$, as required.
\end{proof}

We are now ready to prove the main theorem of this subsection.

\begin{theorem} Let $S$ and $T$ be semigroups with local units.
Then if $S$ and $T$ are Morita equivalent then their Cauchy completions are equivalent.
\end{theorem}
\begin{proof} Let $G$ be the functor of the equivalence that maps  $S-\mbox{\bf FAct}$ to $T-\mbox{\bf FAct}$
and let $H$ be its companion functor going in the opposite direction.
If $M$ is an indecomposable projective in  $S-\mbox{\bf FAct}$ then $G(M)$ is an indecomposable projective in $T-\mbox{\bf FAct}$.
Thus $G$ maps the full subcategory of indecomposable projectives in $S-\mbox{\bf FAct}$
to the full subcategory of indecomposable projectives in $T-\mbox{\bf FAct}$,
and $H$ does the same in the opposite direction.
Thus the full subcategory of indecomposable projectives in $S-\mbox{\bf FAct}$
is equivalent to 
the full subcategory of indecomposable projectives in $T-\mbox{\bf FAct}$.
By Proposition~3.3, every indecomposable projective in  $S-\mbox{\bf FAct}$
is isomorphic to one of the form $Se$ for some idempotent $e$.
Let $\mathcal{IP}_{S}$ be the full subcategory of $S-\mbox{\bf FAct}$
whose objects are all the left closed $S$-sets of the form $Se$ where $e$ ranges over all idempotents of $S$.
Then the full subcategory of indecomposable projectives in $S-\mbox{\bf FAct}$ is equivalent to $\mathcal{IP}_{S}$. 
Similarly, the full subcategory of indecomposable projectives in $T-\mbox{\bf FAct}$ is equivalent to $\mathcal{IP}_{T}$. 
It follows that $\mathcal{IP}_{S}$ is equivalent to  $\mathcal{IP}_{T}$.

Let $\alpha \colon \: Se \rightarrow Sf$ be a left $S$-homomorphism.
Put $a = (e)\alpha$.
Then $a = (e)\alpha = (ee)\alpha = e(e)\alpha = ea$
and $a \in Sf$ and so $af = a$.
It follows that $a = eaf$.
Also $(r)\alpha = (re)\alpha = r(e)\alpha = ra$.
Thus $\alpha = \rho_{a}$ where $a = (e)\alpha$.
Conversely, if $b = ebf$ then $\rho_{b} \colon \: Se \rightarrow Sf$ is a left $S$-homomorphism.
Now let 
$$Se \stackrel{\alpha}{\rightarrow} Sf \stackrel{\beta}{\rightarrow} Sg$$
be a composable
sequence of left $S$-homomorphisms.
Put $a = \alpha (e)$ and $b = \beta (f)$.
Then $\alpha \beta = \rho_{a}\rho_{b} = \rho_{ab}$.
Define a map from $C(S)$ to $\mathcal{IP}_{S}$
by $(e,a,f)$ maps to $\rho_{a} \colon \: Se \rightarrow Sf$.
Then this defines a functor which is full and faithful and every object in $\mathcal{IP}_{S}$
is actually in the image of the map.
It follows that $C(S)$ is equivalent to $\mathcal{IP}_{S}$
and thus $C(S)$ and $C(T)$ are equivalent.
\end{proof}

%%%%%%%%%%%%%%%%%%%%%%%%%%%%%%%%%%%%%%%%%%%%%%%%%%%%%%%%%%%%%%%%%%%%%%%%%%%%%%%%%%%%%%%%%%%%%%%%%%%%%%%%%%%%%%%%%%%%%%%%%%%%%%%%%%%%%%%%%%%%%%%%%%%%%%%%%%%
\subsection{From  Cauchy completions to enlargements}

The result in this section is the linchpin of the whole theorem.
The method we use is based on an argument of McAlister \cite{M3} which was formalized in \cite{K1}.
This formalization was then refined using \cite{P1}.

\begin{lemma} Let $C$ be a strongly connected category and let $p$ be a consolidation on $C$.
Then $C^{p}$ is a semigroup with local units.
In addition, if $C$ is regular then $C^{p}$ is regular.
\end{lemma}
\begin{proof}
Let $x \in C$ be an arrow from $e$ to $f$.
Then $x \circ e = xp_{e,e}e = xe = x$.
Similarly, $f \circ x = x$.
Thus $C^{p}$ is a semigroup with local units.
Suppose now that $C$ is regular.
Given $x$ an arrow from $e$ to $f$ there is an arrow $x'$ from $f$ to $e$
such that $x = xx'x$ and $x' = x'xx'$.
But $x \circ x' \circ x = xp_{e,e}x'p_{f,f}x = xx'x = x$.
Similarly $x' = x' \circ x \circ x'$.
Thus $C$ a regular category implies $C^{p}$ a regular semigroup.
\end{proof}

Our next definition is a version of our definition of a bipartite category given in \cite{K1}
sharpened up in the light of the notion of `bridge' discussed in \cite{P1}.
Let $C$ be a category.
We say that $C = [A,B]$ is {\em bipartite (with left part $A$ and right part $B$)} if it satisfies the following conditions:
\begin{description}

\item[{\rm (B1)}] $C$ has full disjoint subcategories $A$ and $B$ such that $C_{o} = A_{o} \cup B_{o}$.

\item[{\rm (B2)}] For each identity $e \in A_{o}$ there exists an isomorphism $x$ with domain $e$ and codomain in $B_{o}$;
for each identity $f \in B_{o}$ there exists an isomorphism $y$ with domain $f$ and codomain in $A_{o}$.

\end{description}

The category $C$ is a disjoint union of four kinds of arrows:
those in $A$; those in $B$; those starting in $A_{o}$ and ending in $B_{o}$;
and those starting in $B_{o}$ and ending in $A_{o}$. 
On this basis, each arrow of $C$ can be assigned one of four {\em types}: $AA,BB,BA,AB$, respectively.
These types multiply as a rectangular band: the type of a product is the product of the types and so is determined by the first and last element of the product.
Observe that if $A$ and $B$ are strongly connected then so too is $C$.
We shall always assume in what follows that $A$ and $B$ are strongly connected.
The following is Theorem~2.2 of \cite{P1}.

\begin{proposition} The categories $A$ and $B$ are equivalent if and only if there is a bipartite category
with left part $A$ and right part $B$.
\end{proposition}

A consolidation $r$ of a bipartite category $C$ induces consolidations on the full subcategories $A$ and $B$.
Thus $A^{r}$ and $B^{r}$ are subsemigroups of $C^{r}$.

\begin{lemma} Let $C = [A,B]$ be a bipartite category and let $r$ be a consolidation defined on $C$.
Then $C^{r}$ is an enlargement of both $A^{r}$ and $B^{r}$.
\end{lemma}
\begin{proof}
We shall prove that $C^{r}$ is an enlargement of $A^{r}$.
The proof that $C^{r}$ is an enlargement of $B^{r}$ follows by symmetry.
We have to prove that
$A^{r} \circ C^{r} \circ A^{r} = A^{r}$ 
and
$C^{r} = C^{r} \circ A^{r} \circ C^{r}$.

Observe that $A^{r}$ has local units and so 
$A^{r} \subseteq A^{r} \circ C^{r} \circ A^{r}$ is immediate.
We prove the reverse inclusion.
Let $a,a' \in A$ and $c \in C$.
Then $a \circ c \circ a' = ar_{e,f}c_{e',f'}ra'$ for suitable identities $e,e',f,f'$.
But the element on the righthand side begins and ends in $A$,
and $A$ is a full subcategory of $C$ and so belongs to $A$, as required.
Thus we have proved the first equality.

Observe that 
$C^{r} \circ A^{r} \circ C^{r} \subseteq C^{r}$ always.
We prove the reverse inclusion.
Let $c \in C$ be arbitrary.
There are four cases to consider.\\

\noindent
1.~Suppose that $c \in A$.
Then $A^{r}$ has local units and so we can write $c = e \circ c \circ f$ where $e,f$ are identities in $A$ and so in $C$.\\

\noindent
2.~Suppose that $c \in B$. Then there are isomorphisms $\alpha$ and $\beta$ such that
$\alpha c \beta^{-1}$ begins and ends in $A$.
Thus it must belong to $A$ because $A$ is a full subcategory.
Put $a = \alpha c \beta^{-1}$.
Then $c = \alpha^{-1}a\beta = \alpha^{-1} \circ a \circ \beta$.\\

\noindent
3.~Suppose that $c$ begins in $A_{o}$ and ends in $B_{o}$.
Then $c = cee$ where $e$ is the domain of $c$ in $A_{o}$.\\

\noindent
4.~Suppose that $c$ begins in $B_{o}$ and ends in $A_{o}$.
Then $c = eec$ where $e$ is the codomain of $c$ in $A_{o}$.\\

Thus $C^{r} \subseteq C^{r} \circ A^{r} \circ C^{r}$ in all cases.
\end{proof}

\begin{lemma} Let $C$ be a category in which each local monoid has the property that the idempotents generate a regular subsemigroup.
Then the set of regular elements of $C$ forms a subcategory.
\end{lemma}
\begin{proof}
Let $e \stackrel{x}{\leftarrow} f$ and $f \stackrel{y}{\leftarrow} g$ be regular elements
with inverses $x'$ and $y'$ respectively.
The idempotents $x'x$ and $yy'$ belong to the local monoid at $f$.
By assumption, the sandwich set $S(x'x,yy')$ is non-empty; see Proposition~2.5.1 of \cite{H}.
Let $i \in S(x'x,y'y)$.
Recall that by definition,
$i$ is an idempotent and $i$ is an inverse of $(x'x)(y'y)$.
Consider the element $y'ix'$.
Then 
$$xy(y'ix')xy 
= x(yy')i(x'x)y 
= x(x'x)(yy')i(x'x)(yy')y
= x(x'x)(yy')y
=
xy,$$
and so 
$xy$ is regular.
\end{proof}

Isomorphisms are regular and so by Lemma~3.8 we have the following.

\begin{corollary} Let $C = [A,B]$ be a bipartite category.
Suppose that $A$ and $B$ are both regular.
Then we can assume that $C$ is also regular.
\end{corollary}

%%%%%%%%%%%%%%%%%%%%%%%%%%%%%%%%%%%%%%%%%%%%%%%%%%%%%%%%%%%%%%%%%%%%%%%%%%%%%%%%%%%%%%%%%%%%%%%%%%%%%%%%%%%%%%%%%%%%%%%%%%%%%%%%%%%%%%%%%%%%%%%%%%%%%%%%%%%%%%%%%%%
Let $C = [A,B]$ be a bipartite category, let $p$ be a consolidation on $A$, and let $q$ be a consolidation on $B$.
Choose an identity $i_{0} \in A_{o}$ and an isomorphism $\xi$ with domain $i_{0}$ and codomain $j_{0} \in B_{o}$.
Define a consolidation $r$ on $C$ as follows:
$$
r_{e,f} = \left\{
\begin{array}{ll}
p_{e,f} & \mbox{ if $e,f \in A_{o}$}\\
q_{e,f} & \mbox{if $e,f \in B_{o}$}\\
q_{e,j_{0}}\xi p_{i_{0},f} & \mbox{if $e \in B_{o}, f \in A_{o}$}\\
p_{e,i_{0}}\xi^{-1}q_{j_{0},f} & \mbox{ if $e \in A_{o}, f \in B_{o}$}
\end{array}
\right.
$$
In other words, $r$ agrees with $p$ and $q$ on $A$ and $B$ respectively, and then uses $\xi$ to do the simplest possible thing to define it on the whole of $C$.
We say that $r$ is the {\em natural extension of $p$ and $q$ to $C$ via $\xi$.}
The following lemma, where we assume the above setup, is proved by means of routine verifications.
The results are expressed in terms of types. 
Thus $AB \circ AA$ means the product of an element of type $AB$ with an element of type $AA$.

\begin{lemma} \mbox{}
\begin{enumerate}

\item $AB \circ \xi \circ \xi^{-1} \circ AA = AB \circ AA$.

\item $AA \circ \xi \circ \xi^{-1} \circ BA = AA \circ BA$.

\item $AA \circ \xi^{-1} \circ BB = AA \circ BB$.

\item $BB \circ \xi \circ AA = BB \circ AA$.

\end{enumerate}
\end{lemma}

%%%%%%%%%%%%%%%%%%%%%%%%%%%%%%%%%%%%%%%%%%%%%%%%%%%%%%%%%%%%%%%%%%%%%%%%%%%%%%%%%%%%%%%%%%%%%%%%%%%%%%%%%%%%%%%%%%%%%%%%%%%%%%%%%%%%%%%%%%%%%%%%%%%%%%%%%%%%%%%%%%
\begin{proposition} Let $C = [A,B]$ be a bipartite category where $A$ and $B$ are strongly connected.
Let $p$ be a consolidation on $A$, and $q$ a consolidation on $B$.
Choose an identity $i_{0} \in A_{o}$ and an isomorphism $\xi$ with domain $i_{0}$ and codomain $j_{0} \in B_{o}$,
and let $r$ be the natural extension of $p$ and $q$ to $C$ via $\xi$.

Let $\pi_{1}$ be a congruence on $A^{p}$,
and let $\pi_{2}$ be a congruence on $B^{q}$.
Let $\pi$ be the congruence on $C^{q}$ generated by $\pi_{1} \cup \pi_{2}$.
\begin{enumerate}

\item $\pi \cap (A \times A) = \pi_{1}$ if the following three conditions hold: 
\begin{description}
\item[{\rm (i)}] $(a,a') \in \pi_{1} \Rightarrow (\xi^{-1} \circ a, \xi^{-1} \circ a') \in \pi_{1}$.

\item[{\rm (ii)}] $(a,a') \in \pi_{1} \Rightarrow (a \circ \xi, a' \circ \xi) \in \pi_{1}$.

\item[{\rm (iii)}] $(b,b') \in \pi_{2} \Rightarrow (\alpha \circ b \circ \beta, \alpha \circ b' \circ \beta) \in \pi_{1}$
for all isomorphisms $\alpha$ and $\beta$ where $\alpha$ is of type $AB$ and $\beta$ is of type $BA$.

\end{description}

\item $\pi \cap (B \times B) = \pi_{2}$ if the following three conditions hold:
\begin{description}
\item[{\rm (i)}] $(b,b') \in \pi_{2} \Rightarrow (\xi \circ b, \xi \circ b') \in \pi_{2}$

\item[{\rm (ii)}] $(b,b') \in \pi_{2} \Rightarrow (b \circ \xi^{-1}, b' \circ \xi^{-1}) \in \pi_{2}$.

\item[{\rm (iii)}] $(a,a') \in \pi_{1} \Rightarrow (\alpha \circ a \circ \beta, \alpha \circ a' \beta) \in \pi_{2}$
for all isomorphisms $\alpha$ and $\beta$ where $\alpha$ is of type $BA$ and $\beta$ is of type $AB$.

\end{description}

\end{enumerate}
\end{proposition}
\begin{proof}
We shall prove (1); the proof of (2) follows by symmetry.

Let $a,a' \in A$ such that $a \, \pi \, a'$.
Then there is a sequence of elementary transitions
$$a = z_{1} \rightarrow z_{2} \rightarrow \ldots \rightarrow z_{n-1} \rightarrow z_{n} = a'$$
where $z_{i} = x_{i} \circ u_{i} \circ y_{i}$, $z_{i+1} = x_{i} \circ v_{i} \circ y_{i}$ and $(u_{i},v_{i}) \in \pi_{1} \cup \pi_{2}$.

We first show that $z_{1},z_{2}, \ldots, z_{n} \in A$.
Recall that each element of $C$ has one of four types and that the type of a product is the (rectangular band) product of the types.
Since $z_{1} = a \in A$ we must have that $x_{1}$ and $y_{1}$
have types $A \ast$ and $\ast A$ where $\ast$ can be either $A$ or $B$.
It follows that $z_{2} \in A$.
We can now repeat this argument to get that all the remaining $z_{i} \in A$, as claimed.

Let $z_{i} \rightarrow z_{i+1}$ be an elementary transition where $z_{i},z_{i+1} \in A$.
We shall prove that in fact $z_{i} \pi_{1} z_{i+1}$, which will establish our claim.
There are two cases to consider.\\

\noindent
{\em Case~1.} We suppose that $u_{i} \, \pi_{1} \, v_{i}$.
Given that $z_{i}, z_{i+1} \in A$, there are four possibilities for the types of $x_{i}$ and $y_{i}$, respectively:
\begin{description}

\item[{\rm (1)}] $AA$ and $AA$.

\item[{\rm (2)}] $AB$ and $AA$.

\item[{\rm (3)}] $AA$ and $BA$.

\item[{\rm (4)}] $AB$ and $BA$.

\end{description}
We shall deal with each of these possibilities in turn.\\

\noindent
(1) We are given that $x_{i}$ and $y_{i}$ are both of type $AA$ and $u_{i} \, \pi_{1} \, v_{i}$.
Thus $z_{i} \pi_{1} z_{i+1}$, as required.\\

\noindent
(2) We are given that $x_{i}$ is of type $AB$ and $y_{i}$ is of type $AA$.
Since $y_{i} \in A$ we have that $u_{i} \circ y_{i} \, \pi_{1} \, v_{i} \circ y_{i}$.
By condition 1(i), we have that
$\xi^{-1} \circ u_{i} \circ y_{i} \, \pi_{1} \, \xi^{-1} \circ v_{i} \circ y_{i}$.
By fullness, $x_{i} \circ \xi \in A$.
Thus
$$x_{i} \circ \xi \circ \xi^{-1} \circ u_{i} \circ y_{i} 
\, \pi_{1} \, 
x_{i} \circ \xi \circ \xi^{-1} \circ v_{i} \circ y_{i}.$$
This simplifies according to Lemma~3.10(1).\\

\noindent
(3) We are given that $x_{i}$ is of type $AA$ and $y_{i}$ is of type $BA$.
Since $x_{i} \in A$ we have that $x_{i} \circ u_{i} \, \pi_{1} \, x_{i} \circ v_{i}$.
By condition 1(ii), we have that
$x_{i} \circ u_{i} \circ \xi \, \pi_{1} \, x_{i} \circ v_{i} \circ \xi$.
By fullness, $\xi^{-1} \circ y_{i} \in A$. 
Thus
$$x_{i} \circ u_{i} \circ \xi \circ \xi^{-1} \circ y_{i} 
\, \pi_{1} \, 
x_{i} \circ v_{i} \circ \xi \circ \xi^{-1} \circ y_{i}.$$
This simplifies by Lemma~3.10(2).\\

\noindent
(4) This case follows by (2) and (3) above.\\

\noindent
{\em Case~2.} We suppose that $u_{i} \, \pi_{2} \, v_{i}$.
Given that $z_{i},z_{i+1} \in A$, there are four possibilities for the types of $x_{i}$ and $y_{i}$, respectively:
\begin{description}

\item[{\rm (I)}] $AA$ and $AA$.

\item[{\rm (II)}] $AB$ and $AA$.

\item[{\rm (III)}] $AA$ and $BA$.

\item[{\rm (IV)}] $AB$ and $BA$.

\end{description}
We shall deal with each of these possibilities in turn.\\

\noindent
(I) Let $x_{i}$ and $y_{i}$ both have type $AA$.
By condition 1(iii), we have that
$\xi^{-1} \circ u_{i} \circ \xi \, \pi_{1} \, \xi^{-1} \circ v_{i} \circ \xi$.
But $x_{i},y_{i} \in A$ and so
$$x_{i} \circ \xi^{-1} \circ u_{i} \circ \xi \circ y_{i} \, \pi_{1} \, x_{i} \circ \xi^{-1} \circ v_{i} \circ \xi \circ y_{i}.$$
This simplifies according to Lemma~3.10(3) and (4).\\

\noindent
(II) Let $x_{i}$ have type $AB$ and $y_{i}$ have type $AA$.
Let the domain of $x_{i}$ be $e$.
By assumption, there is an isomorphism $\alpha$ that starts at $e$ and ends in $A_{o}$.
By condition 1(iii), we have that
$\alpha \circ u_{i} \circ \xi \, \pi_{1} \, \alpha \circ v_{i} \circ \xi$.
Now $y_{i} \in A$ and so 
$\alpha \circ u_{i} \circ \xi \circ y_{i} \, \pi_{1} \, \alpha \circ v_{i} \circ \xi \circ y_{i}$.
Thus by Lemma~3.10(4), we have that
$\alpha \circ u_{i} \circ y_{i} \, \pi_{1} \, \alpha \circ v_{i} \circ y_{i}$.
Now $x_{i} \circ \alpha^{-1} \in A$ by fullness.
Thus
$$
x_{i} \circ \alpha^{-1} \circ
\alpha \circ u_{i} \circ y_{i} 
\, \pi_{1} \, 
x_{i} \circ \alpha^{-1} \circ
\alpha \circ v_{i} \circ y_{i}.
$$
But $x_{i} \circ \alpha^{-1} \circ \alpha = x_{i} \circ e = x_{i}$,
and so we get the required result.\\

\noindent
(III) Let $x_{i}$ have type $AA$ and let $y_{i}$ have type $BA$.
Let the codomain of $y_{i}$ be $e$.
By assumption, there is an isomorphism $\beta$ that ends at $e$ and starts in $A_{o}$.
By condition 1(iii), we have that
$\xi^{-1} \circ u_{i} \circ \beta \, \pi_{1} \, \xi^{-1} \circ v_{i} \circ \beta$.
Now $x_{i} \in AA$ and so 
$x_{i} \circ \xi^{-1} \circ u_{i} \circ \beta \, \pi_{1} \, x_{i} \circ \xi^{-1} \circ v_{i} \circ \beta$.
Thus by Lemma~3.10(3), we have that
$x_{i} \circ u_{i} \circ \beta \, \pi_{1} \, x_{i} \circ v_{i} \circ \beta$.
Now $\beta^{-1} \circ y_{i} \in A$ by fullness.
Thus
$$x_{i} \circ u_{i} \circ \beta \circ \beta^{-1} \circ y_{i} \, \pi_{1} \, x_{i} \circ v_{i} \circ \beta \circ \beta^{-1} \circ y_{i}.$$
But $\beta \circ \beta^{-1} \circ y_{i} = e \circ y_{i} = y_{i}$,
and so we get the required result.\\

\noindent
(IV) Let $x_{i}$ have type $AB$ and let $y_{i}$ have type $BA$.
Let the domain of $x_{i}$ be $e$ and let the codomain of $y_{i}$ be $f$.
By condition 1(iii), we have that 
$\alpha \circ u_{i} \circ \beta \, \pi_{1} \, \alpha \circ v_{i} \circ \beta$.
By fullness, $x_{i} \circ \alpha^{-1} \in A$ and $\beta^{-1} \circ y_{i} \in A$.
Thus
$$ x_{i} \circ \alpha^{-1} \circ
\alpha \circ u_{i} \circ \beta 
\, \pi_{1} \, 
\alpha \circ v_{i} \circ \beta
\circ \beta^{-1} \circ y_{i}.
$$
But this simplifies as before to the required result.
\end{proof}

Enlargements are preserved under homomorphisms by Proposition~2.9 of \cite{K1}.

\begin{proposition} Let $C = [A,B]$ be a bipartite category where $A$ and $B$ are strongly connected.
Let $p$ be a consolidation on $A$, and $q$ a consolidation on $B$.
Choose an identity $i_{0} \in A_{o}$ and an isomorphism $\xi$ with domain $i_{0}$ and codomain $j_{0} \in B_{o}$,
and let $r$ be the natural extension of $p$ and $q$ to $C$ via $\xi$.
Let $S$ be a homomorphic image of $A^{p}$ by a map with kernel $\pi_{1}$,
and let $T$ be a homomorphic image of $B^{q}$ by a map with kernel $\pi_{2}$.
Let $\pi$ be the congruence on $C^{q}$ generated by $\pi_{1} \cup \pi_{2}$.
Put $R = C^{r}/\pi$.
Suppose that conditions (1) and (2) of Proposition~3.11 hold.
Then $R$ is an enlargement of both $S$ and $T$.
\end{proposition}

%%%%%%%%%%%%%%%%%%%%%%%%%%%%%%%%%%%%%%%%%%%%%%%%%%%%%%%%%%%%%%%%%%%%%%%%%%%%%%%%%%%%%%%%%%%%%%%%%%%%%%%%%%%%%%%%%%%%%%%%%%%%%%%%%%%%%%%%%%%%%%%%%%%%%%%%%%%%%%%%%%
We now apply these results to the problem in hand.

\begin{theorem} Let $S$ and $T$ be semigroups with local units.
If the categories $C(S)$ and $C(T)$ are equivalent then $S$ and $T$ have a joint enlargement
which can be chosen to be regular if $S$ and $T$ are both regular.
\end{theorem}
\begin{proof}
Let $C(S)$ and $C(T)$ be equivalent categories.
By Proposition~3.6, we can find a bipartite category $C = [C(S),C(T)]$.
Both $C(S)$ and $C(T)$ are strongly connected and so $C$ is strongly connected.
We now make the following definitions.

\begin{itemize}

\item The identities of $C(S)$ are of the form $(e,e,e)$ where $e$ is an idempotent of $S$.
We abbreviate $(e,e,e)$ by $\overline{e}$.
On $C(S)$ we define the consolidation $p$ by $p_{\overline{e},\overline{f}} = (e,ef,f)$.
The function $\pi_{1}^{\natural} \colon \: C(S)^{p} \rightarrow S$ given by
$(e,s,f) \mapsto s$ is a surjective homomorphism.

\item The identities of $C(T)$ are of the form $(i,i,i)$ where $i$ is an idempotent of $T$.
We abbreviate $(i,i,i)$ by $\mathbf{i}$.
On $C(T)$ we define the consolidation $q$ by $q_{\mathbf{i},\mathbf{j}} = (i,ij,j)$.
The function $\pi_{2}^{\natural} \colon \: C(T)^{q} \rightarrow T$ given by
$(i,t,j) \mapsto t$ is a surjective homomorphism.

\end{itemize}

Let $\overline{e_{0}}$ be any identity in $C(S)$.
Since $C$ is bipartite, there is an isomorphism $\xi \in C$ with domain $\overline{e_{0}}$
and codomain $\mathbf{f_{0}}$ for some identity in $C(T)$.
Let $r$ be a natural extension of $p$ and $q$ to $C$ defined using this $\xi$.
We now verify that the conditions of Proposition~3.11(1) hold;
that those of (2) also hold follows by symmetry.

Condition (i). Let $(e,s,f) \, \pi_{1} \,(e',s,f')$.
Then simple calculations show that
$\xi^{-1} \circ (e,s,f) = (e_{0},e_{0}s,f)$ 
and
$\xi^{-1} \circ (e',s,f') = (e_{0},e_{0}s,f')$. 
Hence $\xi^{-1} \circ (e,s,f) \, \pi_{1} \,\xi^{-1} \circ (e',s,f')$.

Condition (ii). Let $(e,s,f) \, \pi_{1} \,(e',s,f')$.
Then simple calculations show that
$(e,s,f) \circ \xi = (e,se_{0},e_{0})$ 
and
$(e',s,f') \circ \xi = (e,se_{0},e_{0})$. 
Hence $(e,s,f) \circ \xi \, \pi_{1} \, (e',s,f') \circ \xi$.

Condition (iii). Let $(i,t,j) \, \pi_{2} \, (i',t,j')$.
Let $\overline{f} \stackrel{\alpha}{\rightarrow} \mathbf{e}$ 
and 
$\mathbf{e'} \stackrel{\beta}{\rightarrow} \overline{f'}$
be isomorphisms in $C$.
Then simple calculations show that
$\alpha \circ (i,t,j) \circ \beta = \alpha (f,ftf',f')\beta$
and
$\alpha \circ (i',t,j') \circ \beta = \alpha (f,ftf',f')\beta$.
Thus these two elements are actually equal and so clearly $\pi_{1}$-related.

By Proposition~3.12,
$R = C^{r}/\pi$ is a semigroup with local units that is an enlargement of (isomorphic copies of) $S$ and $T$.

If both $C(S)$ and $C(T)$ were regular, 
then we could assume that $C$ was regular and so $C^{r}$ was regular by Lemma~3.9.
Hence $R$ would be regular, as required.
\end{proof}

%%%%%%%%%%%%%%%%%%%%%%%%%%%%%%%%%%%%%%%%%%%%%%%%%%%%%%%%%%%%%%%%%%%%%%%%%%%%%%%%%%%%%%%%%%%%%%%%%%%%%%%%%%%%%%%%%%%%%%%%%%%%%%%%%%%%%%%
\subsection{From  enlargements to Morita contexts}

\begin{proposition} Let $S$ and $T$ have a joint enlargement $R$.
Then one can construct a unitary Morita context with surjective maps $(S,T,P,Q,\langle -,-\rangle,[-,-])$.
\end{proposition}
\begin{proof}
Put $P = SRT$ and $Q = TRS$.
Then under left and right multiplication, $P$ is an $(S,T)$-biact and $Q$ is a $(T,S)$-biact.
The fact that the left $S$-act $P$ is unitary follows from the fact that $SP = SSRT = SRT = P$.
We prove that $P$ is a closed left $S$-act.
It is enough to prove that
$$\mu \colon \: S \otimes P \rightarrow P$$
is injective.
Let $s \otimes s_{1}r_{1}t_{1}, s' \otimes s_{2}r_{2}t_{2} \in S \otimes SR$ be such that $ss_{1}r_{1}t_{1} = s's_{2}r_{2}t_{2}$.
We prove that $s \otimes s_{1}r_{1}t_{1} = s' \otimes s_{2}r_{2}t_{2}$.
Let $f \in T$ be an idempotent such that $t_{1}f = t_{1}$.
Then 
$$ss_{1}r_{1}t_{1}f = ss_{1}r_{1}t_{1} \mbox{ and }ss_{2}r_{2}t_{2}f = ss_{2}r_{2}t_{2}.$$
Since $f \in T \subseteq R = RSR$ we can write
$$f = r_{3}s_{3}r_{4}$$
where $r_{3},r_{4} \in R$ and $s_{3} \in S$.
Since $s_{3},s' \in S$ there exist idempotents $e,i \in S$ such that
$$s_{3}e = s_{3} \mbox{ and } is' = s'.$$
We now calculate
\begin{eqnarray*}
s \otimes s_{1}r_{1}t_{1} &=& s \otimes s_{1}r_{1}t_{1}f\\
                    &=& s \otimes s_{1}r_{1}t_{1}(r_{3}s_{3}r_{4})f\\
                           &=& s \otimes (s_{1}r_{1}t_{1}r_{3}s_{3})(er_{4}f)\\
                           &=& s(s_{1}r_{1}t_{1}r_{3}s_{3}) \otimes er_{4}f \\
                           &=& (ss_{1}r_{1}t_{1})(r_{3}s_{3}) \otimes er_{4}f\\
                           &=& (s's_{2}r_{2}t_{2})(r_{3}s_{3}) \otimes er_{4}f\\
                           &=& s'(s_{2}r_{2}t_{2}r_{3}s_{3}) \otimes er_{4}f\\
                           &=& s' \otimes (s_{2}r_{2}t_{2}r_{3}s_{3})er_{4}f\\
                           &=& s' \otimes (s_{2}r_{2}t_{2})(r_{3}s_{3}r_{4})f\\
                           &=& s' \otimes (s_{2}r_{2}t_{2})f\\
                           &=& i \otimes (s's_{2}r_{2}t_{2})f\\
                           &=& i \otimes s's_{2}r_{2}t_{2}\\
                           &=& s' \otimes s_{2}r_{2}t_{2}
\end{eqnarray*}
as required.
It follows that $P$ is a closed left $S$-act.
The map
$$\langle -,- \rangle \colon \:P \otimes Q \rightarrow S$$
is defined by $\langle p,q \rangle = pq$ where $p \in SRT$ and $q \in TRS$.
Observe that $pq \in (SRT)(TRS) = S(RTTR)S \subseteq SRS = S$ and is well-defined.
It is clearly an $(S,S)$-homomorphism.
The fact that this map is surjective follows from the fact that
$$S = SRS = S(RTR)S = (SRT)(TRS) = PQ.$$
The map 
$$[-,-] \colon \: Q \otimes P \rightarrow T$$
is defined by $[q,p] = qp$.
This is clearly a $(T,T)$-homomorphism
and surjective by a similar argument to the above.
It is now immediate that we have defined a unitary Morita context with surjective maps.
\end{proof}

%%%%%%%%%%%%%%%%%%%%%%%%%%%%%%%%%%%%%%%%%%%%%%%%%%%%%%%%%%%%%%%%%%%%%%%%%%%%%%%%%%%%%%%%%%%%%%%%%%%%%%%%%%%%%%%%%%%%%%%%%%
\subsection{From  Morita contexts to Morita equivalence}

The following was first proved as Lemma~8.1 of \cite{T1}.
We give an alternative proof.

\begin{lemma}
Let $(S,T,P,Q,\langle -,- \rangle,[-,-])$ be a unitary Morita context with $\langle -,-\rangle$ and $[-,-]$ surjective. 
Then $\langle -,- \rangle$ and $[-,-]$ are injective.
\end{lemma} 
\begin{proof}
We will prove that if $[-,-]: Q \otimes P \rightarrow T$ is surjective then it is injective.
The result for $\langle -,-\rangle$ can be proved similarly.
Let $q \otimes p, q' \otimes p' \in Q \otimes P$ such that $[q,p] = [q',p']$.
Let $e,f \in E(T)$ such that $eq = q$ and $p'f = p'$. 
Since $[,]$ is surjective,
there are $e_{1} \otimes e_{2},f_{1} \otimes f_{2} \in Q \otimes P$ such
that $[e_{1},e_{2}] = e$ and $[f_{1},f_{2}] = f$. 
We have that
$$q \otimes p = (eq) \otimes p = ([e_{1},e_{2}]q) \otimes p.$$
But $[e_{1},e_{2}]q = e_{1}\langle e_{2},q \rangle$
and so
$$([e_{1},e_{2}]q) \otimes p = (e_{1}\langle e_{2},q \rangle) \otimes p = e_{1} \otimes \langle e_{2},q \rangle p.$$
But
$$\langle e_{2},q \rangle p = e_{2}[q,p] = e_{2}[q',p']$$
and so
$$ e_{1} \otimes \langle e_{2},q \rangle p = e_{1} \otimes e_{2}[q',p'].$$
However
$$ e_{1} \otimes (e_{2}[q',p']) 
= e_{1} \otimes (\langle e_{2},q' \rangle p')
=e_{1}\langle e_{2},q' \rangle \otimes p'
= ([e_{1},e_{2}]q') \otimes p'
= eq' \otimes p'.$$
We have proved that $q \otimes p = e(q' \otimes p')$,
and we may similarly prove that $q' \otimes p' = (q \otimes p)f$.
Observe that $e(q \otimes p) = q \otimes p$.
Hence
$$q \otimes p = e(q' \otimes p') = e (q \otimes p)f = (q \otimes p)f = q' \otimes p'.$$
It follows that $q \otimes p = q' \otimes p'$, as required.
\end{proof}

It follows from the above result that $P \otimes Q \cong S$, as an $(S,S)$-biact, and $Q \otimes P \cong T$, as a $(T,T)$-biact.
The  following was first proved as Theorem~8.3 of \cite{T1}.

\begin{proposition} Let $(S,T,P,Q,\langle -,-\rangle,[-,-])$ be a unitary Morita context with surjective maps.
Then the categories $S-\mbox{\bf FAct}$ and $T-\mbox{\bf FAct}$ are equivalent via the functors 
$$Q \otimes - \colon S-\mbox{\bf FAct} \rightarrow T-\mbox{\bf FAct}$$
and  
$$P \otimes - \colon T-\mbox{\bf FAct} \rightarrow S-\mbox{\bf FAct}.$$
\end{proposition}
\begin{proof}
Let $M$ be a closed left $S$-act.
Then we may form $Q \otimes M$.
Observe that
$$T \otimes (Q \otimes M) \cong (T \otimes Q) \otimes M \cong Q \otimes M$$
because $Q$ is left closed.
Thus $Q \otimes M$ is left closed.
It follows that we have well-defined functors in each direction.
It remains to show that they form an equivalence of categories.
Let $M$ be a closed left $S$-act.
Then 
$$P \otimes (Q \otimes M) \cong (P \otimes Q) \otimes M \cong S \otimes M \cong M$$ 
using the remark following Lemma~3.15.
We therefore have a left $S$-isomorphism 
$$\alpha_{M} \colon \: P \otimes (Q \otimes M) \rightarrow M$$ 
which maps $p \otimes (q \otimes m)$ to $(\langle p,q \rangle \otimes m)\mu_{M} = \langle p, q \rangle m$.
It follows that $\alpha$ is a natural equivalence with all components isomorphisms.
A similar result in the other direction leads to the required equivalence of categories.
\end{proof}

\begin{lemma} Let $(S,T,P,Q,\langle -,- \rangle,[-,-])$ be a Morita context with $\langle -,- \rangle$ and $[-,-]$
surjective. Then $P$ and $Q$ are also closed on the right.
\end{lemma}
\begin{proof}
We have the following isomorphisms of biacts
$$P \otimes T \cong P \otimes (Q \otimes P) \cong (P \otimes Q) \otimes P \cong S \otimes P \cong P.$$
Thus $P$ is also closed on the right.
We may similarly show that $Q$ is closed on the right.
\end{proof}

The following result is not stated by Talwar but now follows immediately by Lemma~3.17.

\begin{theorem} Let $S$ and $T$ be semigroups with local units.
Then the categories $S-\mbox{\bf FAct}$ and  $T-\mbox{\bf FAct}$
are equivalent 
if and only if
the categories
$\mbox{\bf FAct}-S$ and  $\mbox{\bf FAct}-T$ are equivalent.
\end{theorem}

%%%%%%%%%%%%%%%%%%%%%%%%%%%%%%%%%%%%%%%%%%%%%%%%%%%%%%%%%%%%%%%%%%%%%%%%%%%%%%%%%%%%%%%%%%%%%%%%%%%%%%%%%%%%%%%%%%%%%%%%%%%%%%%%%%%%%%%%%%%%%
\section{Proof of Theorem~1.2}\setcounter{theorem}{0}

First we prove a lemma which extracts the key result from the proof of Proposition~1 of \cite{L3}.

\begin{lemma} Let $T$ be an enlargement of $S$ where $T^{2} = T$ and $S^{2} = S$.
Then each idempotent of $T$ is $\mathcal{D}$-related to an idempotent of $S$.
\end{lemma}
\begin{proof} Let $e \in E(T)$. 
By assumption $T = TST$,
and so we can write $e = usv$ where $u,v \in T$ and $s \in S$.
By assumption $S^{2} = S$, and so we may write $s = ab$ where $a,b \in S$.
Put $x = ua$ and $y = bv$.
Then $e = xy$.
It is easy to check that $ye \in V(ex)$.
Put $f = yex$, an idempotent.
Then 
$$(ex)(ye) = e \mbox{ and }(ye)(ex) = f$$
and so $e\, \mathcal{D} \, f$.
But $f = yex = (bv)e(ua) = b(veu)a \in S(TTT)S = STS = S$, as required.
\end{proof}

It follows from the lemma below that local isomorphisms are precisely surjective local isomorphisms followed by enlargements.

\begin{lemma} \mbox{}

\begin{enumerate}

\item Let $\alpha \colon \: S \rightarrow T$ and $\beta \colon \: T \rightarrow U$ be local isomorphisms.
Then $\beta \alpha \colon \: S \rightarrow U$ is a local isomorphism.

\item Let $T$ be a semigroup with local units and let $S$ be a subsemigroup of $T$
also with local units. Then $T$ is an enlargement of $S$ if and only if the embedding of $S$ in $T$
is a local isomorphism.

\item  Let $\theta \colon \: S \rightarrow T$ be a local isomorphism.
Then $T$ is an enlargement of the image of $\theta$.

\end{enumerate}

\end{lemma}
\begin{proof} The proof of (1) is straightforward.

(2) If $T$ is an enlargement of $S$ then (L1) and (L2) are immediate, and (LI3) is a consequence of Lemma~4.1.
To prove the converse, suppose that the embedding of $S$ in $T$ is a local isomorphism.
We prove that $T$ is an enlargement of $S$.
For each pair of idempotents $e,f \in S$ we have that $eSf = eTf$.
Let $a \in STS$.
Then there exists $e,f \in S$ such that $a = eaf$.
Then $a \in eTf = eSf$ and so $a \in S$.
It follows that $S = STS$.
Now let $b \in T$ and let $i \in T$ be an idempotent such that $bi = b$.
By assumption, there exists $e \in S$ such that $i \, \mathcal{D} \, e$.
Thus there exists $x \in T$ and $x' \in V(x)$ such that $x'x = i$ and $xx' = e$.
Hence $b = bi = bx'x = bx'xx'x = (bx')ex \in TST$ and so $T = TST$.

(3) Put $T' = \theta (S)$.
Then $T'$ is a subsemigroup of $T$ with local units.
Let $a',b' \in T'$ and $c' \in T$.
Let $\theta (a) = a'$ and $\theta (b) = b'$.
Let $ae = a$ and $fb = b$ where $e$ and $f$ are idempotents in $S$.
Then $a'c'b' = a' \theta (e)c'\theta(f)b'$.
By assumption there exists $c \in eSf$ such that $\theta (c) = c'$.
Thus $a'c'b' = \theta (a)\theta (c)\theta (b) = \theta (acb)$.
It follows that $a'c'b' \in T'$.
We have shown that $T'TT' \subseteq T'$ and so it follows that $T' = T'TT'$.
  
Let $t \in T$.
Then $t = te$ for some $e \in E(T)$.
By assumption $e \, \mathcal{D} \, \theta (f)$ for some $f \in E(S)$.
Let $x \in T$ and $x' \in V(x)$ such that $x'x = e$ and $xx' = \theta (f)$.
Then $t = te = tx'xx'x = (tx')\theta (f)x \in TT'T$.
We have shown that $T \subseteq TT'T$ and so $T = TT'T$, as required.
\end{proof}

%%%%%%%%%%%%%%%%%%%%%%%%%%%%%%%%%%%%%%%%%%%%%%%%%%%%%%%%%%%%%%%%%%%%%%%%%%%%%%%%%%%%%%%%%%%%%%%%%%%%%%%%%%%%%%%%%%%%%%%%%%%
We may now prove Theorem~1.2.

\begin{proof}
Suppose first that there is a local isomorphism $\psi \colon \: C(S)^{q} \rightarrow T$.
We prove that there is an equivalence between $C(S)$ and $C(T)$.
We shall first construct a full and faithful functor
$\Psi \colon \: C(S) \rightarrow C(T)$.
We abbreviate identities $(e,e,e)$ in $C(S)$ by $\mathbf{e}$.
The identity $\mathbf{e}$ is an idempotent in $C(S)^{q}$ because for consolidations $q_{e,e} = e$.
Thus $\psi (\mathbf{e})$ is an idempotent in $T$.
Define $\Psi (e,e,e) = (\psi (\mathbf{e}),\psi (\mathbf{e}),\psi (\mathbf{e}))$.
Now define 
$$\Psi (e,s,f) = (\psi (\mathbf{e}), \psi (e,s,f) ,\psi(\mathbf{f})).$$
This is well-defined because $(e,e,e) \circ (e,s,f) \circ (f,f,f) = (e,s,f)$.
It is routine to check that $\Psi$ is a functor and it is full and faithful because $\psi$ is a local isomorphism.

It remains to show that $\Psi$ is essentially surjective.
We claim  that for each idempotent $i$ in $T$, there is an idempotent $f' \in T$
such that $i \, \mathcal{D} \, f'$ and $\psi (f,f,f) = f'$ for some idempotent $f \in S$.
Before we prove the claim, we show that it implies that $\Psi$ is essentially surjective.
Let $(i,i,i)$ be an identity in $C(T)$.
By assumption,  $i \, \mathcal{D} \, f'$ and $\psi (f,f,f) = f'$ for some idempotent $f \in S$.
Now $i \, \mathcal{D} \, f'$ iff there exists $x \in T$ and $x' \in V(x)$ such that $x'x = f'$ and $xx' = i$.
It follows that $(i,x,f') \in C(T)$ is an isomorphism linking the identities $(i,i,i)$ and $(f',f',f')$.
By assumption $\Psi (f,f,f) = (f',f',f')$, and so $\Psi$ is essentially surjective.

We now prove the claim.
Let $i \in E(T)$.
By assumption there is an idempotent $e' \in E(T)$ such that $i \, \mathcal{D} \, e'$ and $e'$ is in the image of $\psi$.
Because idempotents lift, there is an idempotent $(e,s,f) \in C(S)^{q}$ such that $\psi (e,s,f) = e'$.
The fact that $(e,s,f)$ is an idempotent in $C(S)^{q}$ means that $s = sq_{f,e}s$.
Thus $s$ is regular.
Since $s = esf$ there is an inverse $s' \in V(s)$ such that $s' = fs'e$.
By construction $(f,s',e) \in C(S)$ and in $C(S)^{q}$ we have that $(f,s',e) \in V((e,s,f))$.
It follows that $(f,s',e)(e,s,f) = (f,s's,f)$ is an idempotent in $C(S)^{q}$ 
and of course  $(e,s,f) \, \mathcal{D} \, (f,s's,f) = (f,j,f)$.
At this point, we use the fact that idempotents split in $C(S)$, but we give the details.
Observe that $(f,j,j) \in V((j,j,f))$,
and that
$$(f,j,f) = (f,j,j)(j,j,f)
\mbox{ and }
(j,j,j) = (j,j,f)(f,j,j)$$
giving $(f,j,f) \, \mathcal{D} \, (j,j,j)$.
It follows that
$(e,s,f) \, \mathcal{D} \, (j,j,j)$.
Put $\psi (j,j,j) = f'$.
Then
$$\psi (j,j,j) = f' \, \mathcal{D}\, i,$$ as required.

To prove the converse, let $S$ and $T$ have common enlargement $R$.
We shall prove that there is a subsemigroup $T'$ of $T$ such that $T$ is an enlargement of $T'$,
and there is a consolidation $q$ on $C(S)$ and a surjective local isomorphism $\psi \colon \: C(S)^{q} \rightarrow T'$.
It follows then from our results on local isomorphisms in Lemma~4.2 that we therefore have a local isomorphism from $C(S)^{q}$ to $T$.

By Lemma~4.1, for each $e \in E(S)$ there exists $f \in E(T)$ such that $e \, \mathcal{D} \, f$.
Thus there exists $x_{e} \in R$ and $x_{e}' \in V(x_{e})$ such that $x_{e}'x_{e} = e$ and $x_{e}x_{e}' = f$.
Define a consolidation $q$ on $C(S)$ by $q_{i,j} = x_{i}'x_{j}$.
Observe that $q_{i,i} = i$ and $q_{i,j} \in iSj$.
Thus we may form the semigroup $C(S)^{q}$.
Define $\psi \colon \: C(S)^{q} \rightarrow T$ by $\psi (i,a,j) = x_{i}ax_{j}'$.
Observe that
$\psi (i,a,j) \psi (k,b,l) = x_{i}ax_{j}'x_{k}bx_{l}'$ 
whereas
$$\psi ((i,a,j) \circ (k,b,l)) = \psi (i,aq_{j,k}b,l) = x_{i}aq_{j,k}bx_{l}' = x_{i}ax_{j}'x_{k}bx_{l}',$$
and so $\psi$ is a homomorphism.
Let $T'$ be the image of $\psi$.

We prove that idempotents lift along $\psi$.
Let $\psi (i,a,j) = x_{i}ax_{j}' = f$, where $f \in E(T')$.
Then $x_{i}ax_{j}'x_{i}ax_{j}' = x_{i}ax_{j}'$.
Thus $ax_{j}'x_{i}a = a$ and so $(i,a,j)$ is an idempotent in $C(S)^{q}$.
The proof that $\psi$ is a surjective local isomorphism is straightforward.

We claim that $T$ is an enlargement of $T'$.
Let $t \in T$.
Then $t \in eTf$ for some $e,f \in E(T)$.
Choose $y_{f} \in R$ and $y_{f}' \in V(y_{f})$ such that $y_{f}'y_{f} = f$ and $y_{f}y_{f}' = i \in E(S)$.
Choose $y_{e} \in R$ and $y_{e}' \in V(y_{e})$ such that $y_{e}'y_{e} = e$ and $y_{e}y_{e}' = j \in E(S)$.
Put $t' = x_{j}y_{e}ty_{f}'x_{i}'$ and $s = y_{e}ty_{f}'$.
Then $s \in S$ and $y_{e}'x_{j}'t'x_{i}y_{f} = t$.
But $(j,s,i) \in C(S)$ and $\psi (j,s,i) = t'$.
It follows that $t \in TT'T$ and so $T \subseteq TT'T$.
Clearly $TT'T \subseteq T$ giving $T = TT'T$.

Next observe that $T' \subseteq TT'T$ because $T'$ has local units, being the image of a semigroup having local units.
We show that $T'TT' \subseteq T'$.
By definition, 
$$T' = \{x_{i}ax_{j}' \colon \: a \in iSj \mbox{ where } i,j \in E(S) \}.$$
Consider the product
$$t = (x_{i}ax_{j}')t(x_{k}bx_{l}')$$
where $i,j,k,l \in E(S)$ and $a \in iSj$ and $b \in kSl$.
Then $t = x_{i}(iax_{j}'tx_{k}bl)x_{l}'$ where $s = iax_{j}'tx_{k}bl \in S$.
Thus $t = \psi (i,s,l) \in T'$, as required.
We have shown that $T' = TT'T$ and so $T$ is an enlargement of $T'$.
\end{proof}

%%%%%%%%%%%%%%%%%%%%%%%%%%%%%%%%%%%%%%%%%%%%%%%%%%%%%%%%%%%%%%%%%%%%%%%%%%%%%%%%%%%%%%%%%%%%%%%%%%%%%%%%%%%%%%%%%%%%%%%%%%%%%%%%%%%%%%%%%%%%%%%%%%%%%%%%%%%%
\section{Applications }\setcounter{theorem}{0}

In this section, we shall apply our theory to obtain some concrete results about the Morita theory of regular semigroups.
We begin with a list of Morita invariant properties.
These go back to results obtained for enlargements \cite{L1},
and they were known from the Morita framework to Talwar \cite{T1,T2,T3}.

\begin{proposition} Let $S$ and $T$ be semigroups with local units which are Morita equivalent.
\begin{enumerate}

\item Each local submonoid of $S$ is isomorphic to a local submonoid of $T$, and vice-versa.

\item $S$ is regular if and only if $T$ is regular.

\item The cardinalities of the sets of regular $\mathcal{D}$-classes in $S$ and $T$ are the same.

\item The posets of two-sided ideals in $S$ and $T$ are order-isomorphic.

\item The posets of principal two-sided ideals in $S$ and $T$ are order-isomorphic.

\end{enumerate}
\end{proposition}
\begin{proof}

(1) The categories $C(S)$ and $C(T)$ are equivalent.
The local monoid at the identity $(e,e,e)$ in $C(S)$ is isomorphic to the local submonoid $eSe$.
The result now follows.

(2) It is easy to check that the semigroup $S$ is regular if and only if the category $C(S)$ is regular.
If a pair of categories is equivalent then one is regular if and only if the other is regular.
The categories $C(S)$ and $C(T)$ are equivalent and so the result follows. 

(3) If the categories $C(S)$ and $C(T)$ are equivalent then their groupoids of isomorphisms are equivalent
and, in particular, the number of components in their groupoids of isomorphisms is the same.
It remains to show that the number of components in the groupoid of isomorphisms of $C(S)$ is the same
as the number of regular $\mathcal{D}$-classes of $S$.
The element $(e,s,f) \in C(S)$ is an isomorphism if and only if there is $s' \in V(s)$ such that $s's = f$ and $ss' = e$.
Thus the identities $(e,e,e)$ and $(f,f,f)$ are linked by an isomorphism if and only if $e \, \mathcal{D} \,f$.
Finally, a $\mathcal{D}$-class of $S$ is regular if and only if it contains an idempotent \cite{H}.
The result now follows.

(4) One proof of this result generalizes Theorem~3.3(i) of \cite{L1}.
A more direct proof uses the function $\pi \colon C(S) \rightarrow S$ that maps $(e,s,f)$ to $s$.
This sets up an order isomorphism between the ideals of $S$ and the ideals of $C(S)$.
In addition, if $C$ and $D$ are equivalent categories then their lattices of ideals are order isomorphic.
The result now follows by Theorem~1.1(2). 

(5) A direct proof of this result shows that the bijection set up in (iv) above restricts to a bijection between the posets of principal ideals.
However, we can deduce it simply by applying a lattice-theoretic result.
In the lattice of ideals of a semigroup, the completely join irreducible
elements are the principal ideals, 
and an order isomorphism between lattices
maps completely join irreducible elements bijectively to completely join irreducible elements. 
\end{proof}

The following result was known to Talwar \cite{T1}.
It shows how the theory simplifies radically when at least one of the semigroups is a monoid.

\begin{proposition} Let $S$ be a monoid and $T$ a semigroup with local units.
Then $S$ and $T$ are Morita equivalent if and only if there is an idempotent $f$ in $T$ such that
$T = TfT$ and $fTf$ is isomorphic to $S$.
Thus $T$ is an enlargement of $S$.
\end{proposition}
\begin{proof} Only one direction needs proving. 
Suppose that $S$ and $T$ are Morita equivalent.
By Theorem~1.1 there is a semigroup with local units $R$ which contains $S$ and $T$ as subsemigroups and is an enlargement of both of them.
Let the identity of $S$ be $e$.
By Lemma~4.1, there exists an idempotent $f \in T$ such that $e \, \mathcal{D} \, f$.
Thus there are elements $x,y \in R$ such that $xy = e$ and $yx = f$.
Observe that $eRe \subseteq S = eSe \subseteq eRe$.
Thus $S = eRe$.
Also $fRf \subseteq fTf \subseteq fRf$.
Thus $fTf = fRf$.
But $e \, \mathcal{D} \, f$ implies that $eRe$ and $fRf$ are isomorphic
and so $S$ is isomorphic to $fTf$ \cite{L3}.
Finally, we show that $T = TfT$.
It is enough to show that $T \subseteq TfT$.
Let $t \in T$.
Now $R = RSR$ and so we may write $t = rsr'$ where $r,r' \in R$ and $s \in S$.
Let $i,j \in E(T)$ such that $t = itj$.
Then $t = irsr'j$.
But $s \in S$ and so $s = es = ees$.
Thus 
$$t = irx(yx)ysr'j = (irxf)f(fysr'j) \in (TRRT)f(TRSRT) \subseteq TfT,$$
as required. 
\end{proof}

%%%%%%%%%%%%%%%%%%%%%%%%%%%%%%%%%%%%%%%%%%%%%%%%%%%%%%%%%%%%%%%%%%%%%%%%%%%%%%%%%%%%%%%%%%%%%%%%%%%%%%%%%%%%%%%%%%%%%%%%%%%%%%%%%%%%%%%%%%%%%%%%%%%%%%%
The Morita theory of unital rings provides a framework for understanding the Wedderburn-Artin Theorem \cite{Lam}.
The semigroup analogue of simple artinian rings is the class of completely simple semigroups whose structure was described in the famous Rees-Suschkewitsch Theorem \cite{R}.
Our first theorem, which is well-known, sets the scene for this section by giving a number of equivalent characterizations of completely simple semigroups.
Recall that a semigroup $S$ is said to have a property {\em locally} if each local submonoid $eSe$ has that property.
By the {\em local structure} of a semigroup $S$, 
we mean the structure of the local submonoids $eSe$ as $e$ varies over the set of idempotents of $S$.

\begin{theorem} Let $S$ be a semigroup with local units.
Then the following are equivalent.
\begin{enumerate}

\item $S$ is completely simple.

\item $S$ is regular and locally a group.

\item There is an idempotent $e$ such that $S = SeS$ and $eSe$ is a group.

\item $S$ is Morita equivalent to a group.

\end{enumerate}
\end{theorem}
\begin{proof}

(1)$\Rightarrow$(2). A completely simple semigroup is a simple semigroup with a primitive idempotent.
It is easy to deduce that it must be bisimple and so, since $S$ contains an idempotent, it is regular. 
If $e$ is a primitive idempotent then $eSe$ is a group.
But all idempotents are $\mathcal{D}$-related and so all local submonoids are isomorphic.
It follows that $S$ is regular and locally a group.
For all unproved statements see \cite{H}.

(2)$\Rightarrow$(1). It is easy to show that a semigroup which is regular and locally a group must be simple.
Thus $S$ is regular with a primitive idempotent and so it is completely simple.  

(1)$\Rightarrow$(3). A completely simple semigroup is simple and so $S = SeS$.
We have already proved that all local submonoids are groups.

(3)$\Rightarrow$(1). This says precisely that $S$ is an enlargement of the group $eSe$.
Thus we quickly deduce that $S$ is a simple regular semigroup and we are given that $e$ is a primitive idempotent.

(1)$\Rightarrow$(4). If $S$ is completely simple we have seen that it is an enlargement of a group and so
it is Morita equivalent to a group.

(4)$\Rightarrow$(1). Let $S$ be a Morita equivalent to a group $G$.
Groups are regular and so by Proposition~5.1(2), the semigroup $S$ is regular.
Groups are bisimple and so by Proposition~5.1(3), it follows that $S$ is a bisimple.
By Proposition~5.1(1), each local submonoid of $S$ is isomorphic to a local submonoid of $G$.
But the local submonoids of $G$ are just $G$ itself.
It follows that each local submonoid of $S$ is isomorphic to $G$.
Thus $S$ is a bisimple regular semigroup which is locally a group.
It follows that $S$ is completely simple.
\end{proof}

\begin{remark}
{\em Let $S$ be {\em any} semigroup such that $S = SeS$ and $eSe$ is a group.
Then in fact $S$ is completely simple.
It is enough to show that $S$ is regular.
Let $s \in S$.
Then we can write $s = s_{1}es_{2}$.
Put $a = s_{1}e$ and $b = es_{2}$.
Then $ba \in eSe$, a group.
Thus there is an element $g \in eSe$ such that $gba = bag = e$.
We have that $eS = bagS \subseteq bS = ebS \subseteq eS$.
Thus $e \, \mathcal{R} \, b$.
Similarly $e \, \mathcal{L} \, a$.
But $a = ae \, \mathcal{R} \, ab = s$,
and so $s \, \mathcal{D} \, e$ which implies that $s$ is regular.}
\end{remark}

In a series of papers \cite{M1,M1',M2,Mc,M3}, McAlister set about generalising the theory of completely simple semigroups
in their guise as regular semigroups which are locally groups.
In \cite{M1,M1'}, he concentrated on the locally inverse semigroups.
These are natural generalizations of both completely simple semigroups and inverse semigroups.
In the papers \cite{M2,Mc,M3} he generalized his theory to other classes of regular semigroups described by the structure of their local submonoids.
The following is our interpretation of McAlister's results in \cite{M3}.

\begin{theorem} Let $S$ be a semigroup with local units.
\begin{enumerate}

\item $S$ is Morita equivalent to a group if and only if it is completely simple.

\item $S$ is Morita equivalent to an inverse semigroup if and only if it is regular and locally inverse.

\item $S$ is Morita equivalent to a semilattice if and only if it is regular, locally a semilattice, and $S/\mathcal{J}$ is a meet
semilattice under subset inclusion.

\item $S$ is Morita equivalent to an orthodox semigroup if and only if it is regular and locally orthodox.

\item $S$ is Morita equivalent to an $\mathcal{L}$-unipotent semigroup if and only if it is regular and locally $\mathcal{L}$-unipotent.

\item $S$ is Morita equivalent to an $E$-solid semigroup if and only if it is regular and locally $E$-solid.

\item $S$ is Morita equivalent to a union of groups if and only if it is regular, locally a union of groups, 
and  $S/\mathcal{J}$ is a meet semilattice under subset inclusion.

\end{enumerate}
\end{theorem}
\begin{proof}
(1) was proved as Theorem~5.3.
I shall prove (2) and (3).
The remaining results are proved similarly.\\

\noindent
Proof of (2).
Let $S$ be a semigroup with local units that is Morita equivalent to an inverse semigroup.
By Proposition~5.1(1),(2), we quickly deduce that $S$ is regular and locally inverse.
Conversely, let $S$ be a regular locally inverse semigroup.
In \cite{M1'}, McAlister shows how to construct a consolidation $q$ on $C(S)$ such that $C(S)^{q}$ is an orthodox locally inverse semigroup.
Such a semigroup has a minimum inverse congruence whose associated surjective homomorphism
$\theta \colon \: C(S)^{q} \rightarrow T$ is a local isomorphism to the inverse semigroup $T$.
It follows by Theorem~1.2 that $S$ is Morita equivalent to $T$, an inverse semigroup.\\

\noindent
Proof of (3). 
Let $S$ be a semigroup with local units that is Morita equivalent to a semilattice.
By Proposition~5.1(1),(2),(5) we quickly deduce that $S$ is regular, locally a semilattice and that $S/\mathcal{J}$ is a meet semilattice under subset inclusion.
Conversely, let $S$ be a regular semigroup which is locally a semilattice and for which $S/\mathcal{J}$ is a meet
semilattice under subset inclusion.
Then by Theorem~3.3(ii) of \cite{M1'}, each $eSf$ contains a maximum element.
The consolidation $q_{e,f}$ is defined to be this maximum element.
Let $(e,a,f) \in C(S)^{q}$.
Then $a \in eSf$ and so there is an inverse $a' \in V(a) \cap fSe$.
Thus $(f,a',e) \in C(S)^{q}$.
By construction $a \leq q_{e,f}$ and $a' \leq q_{f,e}$.
We have that $a = aa'a \leq aq_{f,e}a$.
But then $a = aq_{f,e}ai$ for some idempotent $i$.
However $ai = a$ and so in fact $a = aq_{f,e}a$.
We have shown that every element of $C(S)^{q}$ is an idempotent.
It is clearly locally inverse so it is a normal band.
As before, such a semigroup has a minimum inverse congruence whose associated surjective homomorphism
$\theta \colon \: C(S)^{q} \rightarrow T$ is a local isomorphism to the semilattice $T$.
It follows by Theorem~1.2 that $S$ is Morita equivalent to $T$, a semilattice.
\end{proof}

We see that Proposition~5.1 can be used to help find necessary conditions for a semigroup $S$ to be Morita equivalent to a semigroup $T$,
whereas Theorem~1.2 can be used to find sufficient conditions.

One aspect of McAlister's theory we have not touched on is his work on Rees matrix covers.
This, however, is a consequence of Theorem~1.1(3) and the work described in \cite{L3}.

%%%%%%%%%%%%%%%%%%%%%%%%%%%%%%%%%%%%%%%%%%%%%%%%%%%%%%%%%%%%%%%%%%%%%%%%%%%%%%%%%%%%%%%%%%%%%%%%%%%%%%%%%%%%%%%%%%%%%%%%%%%%%%%%%%%%%%%%

\end{document}